\documentstyle[11pt]{article}

\catcode`\@=11

\newfam\msxfam
\newfam\msyfam

\def\hexnumber@#1{\ifcase#1 0\or1\or2\or3\or4\or5\or6\or7\or8\or9\or
	A\or B\or C\or D\or E\or F\fi }

\font\teneuf=eufm10
\font\seveneuf=eufm7
\font\fiveeuf=eufm5
\newfam\euffam
\textfont\euffam=\teneuf
\scriptfont\euffam=\seveneuf
\scriptscriptfont\euffam=\fiveeuf
\def\frak{\relaxnext@\ifmmode\let\next\frak@\else
 \def\next{\Err@{Use \string\frak\space only in math mode}}\fi\next}
\def\goth{\relaxnext@\ifmmode\let\next\frak@\else
 \def\next{\Err@{Use \string\goth\space only in math mode}}\fi\next}
\def\frak@#1{{\frak@@{#1}}}
\def\frak@@#1{\noaccents@\fam\euffam#1}

\edef\msx@{\hexnumber@\msxfam}
\edef\msy@{\hexnumber@\msyfam}

\mathchardef\boxdot="2\msx@00
\mathchardef\boxplus="2\msx@01
\mathchardef\boxtimes="2\msx@02
\mathchardef\square="0\msx@03
\mathchardef\blacksquare="0\msx@04
\mathchardef\centerdot="2\msx@05
\mathchardef\lozenge="0\msx@06
\mathchardef\blacklozenge="0\msx@07
\mathchardef\circlearrowright="3\msx@08
\mathchardef\circlearrowleft="3\msx@09
\mathchardef\rightleftharpoons="3\msx@0A
\mathchardef\leftrightharpoons="3\msx@0B
\mathchardef\boxminus="2\msx@0C
\mathchardef\Vdash="3\msx@0D
\mathchardef\Vvdash="3\msx@0E
\mathchardef\vDash="3\msx@0F
\mathchardef\twoheadrightarrow="3\msx@10
\mathchardef\twoheadleftarrow="3\msx@11
\mathchardef\leftleftarrows="3\msx@12
\mathchardef\rightrightarrows="3\msx@13
\mathchardef\upuparrows="3\msx@14
\mathchardef\downdownarrows="3\msx@15
\mathchardef\upharpoonright="3\msx@16

\mathchardef\downharpoonright="3\msx@17
\mathchardef\upharpoonleft="3\msx@18
\mathchardef\downharpoonleft="3\msx@19
\mathchardef\rightarrowtail="3\msx@1A
\mathchardef\leftarrowtail="3\msx@1B
\mathchardef\leftrightarrows="3\msx@1C
\mathchardef\rightleftarrows="3\msx@1D
\mathchardef\Lsh="3\msx@1E
\mathchardef\Rsh="3\msx@1F
\mathchardef\rightsquigarrow="3\msx@20
\mathchardef\leftrightsquigarrow="3\msx@21
\mathchardef\looparrowleft="3\msx@22
\mathchardef\looparrowright="3\msx@23
\mathchardef\circeq="3\msx@24
\mathchardef\succsim="3\msx@25
\mathchardef\gtrsim="3\msx@26
\mathchardef\gtrapprox="3\msx@27
\mathchardef\multimap="3\msx@28
\mathchardef\therefore="3\msx@29
\mathchardef\because="3\msx@2A
\mathchardef\doteqdot="3\msx@2B

\mathchardef\triangleq="3\msx@2C
\mathchardef\precsim="3\msx@2D
\mathchardef\lesssim="3\msx@2E
\mathchardef\lessapprox="3\msx@2F
\mathchardef\eqslantless="3\msx@30
\mathchardef\eqslantgtr="3\msx@31
\mathchardef\curlyeqprec="3\msx@32
\mathchardef\curlyeqsucc="3\msx@33
\mathchardef\preccurlyeq="3\msx@34
\mathchardef\leqq="3\msx@35
\mathchardef\leqslant="3\msx@36
\mathchardef\lessgtr="3\msx@37
\mathchardef\backprime="0\msx@38
\mathchardef\risingdotseq="3\msx@3A
\mathchardef\fallingdotseq="3\msx@3B
\mathchardef\succcurlyeq="3\msx@3C
\mathchardef\geqq="3\msx@3D
\mathchardef\geqslant="3\msx@3E
\mathchardef\gtrless="3\msx@3F
\mathchardef\sqsubset="3\msx@40
\mathchardef\sqsupset="3\msx@41
\mathchardef\vartriangleright="3\msx@42
\mathchardef\vartriangleleft="3\msx@43
\mathchardef\trianglerighteq="3\msx@44
\mathchardef\trianglelefteq="3\msx@45
\mathchardef\bigstar="0\msx@46
\mathchardef\between="3\msx@47
\mathchardef\blacktriangledown="0\msx@48
\mathchardef\blacktriangleright="3\msx@49
\mathchardef\blacktriangleleft="3\msx@4A
\mathchardef\vartriangle="0\msx@4D
\mathchardef\blacktriangle="0\msx@4E
\mathchardef\triangledown="0\msx@4F
\mathchardef\eqcirc="3\msx@50
\mathchardef\lesseqgtr="3\msx@51
\mathchardef\gtreqless="3\msx@52
\mathchardef\lesseqqgtr="3\msx@53
\mathchardef\gtreqqless="3\msx@54
\mathchardef\Rrightarrow="3\msx@56
\mathchardef\Lleftarrow="3\msx@57
\mathchardef\veebar="2\msx@59
\mathchardef\barwedge="2\msx@5A
\mathchardef\doublebarwedge="2\msx@5B
\mathchardef\angle="0\msx@5C
\mathchardef\measuredangle="0\msx@5D
\mathchardef\sphericalangle="0\msx@5E
\mathchardef\varpropto="3\msx@5F
\mathchardef\smallsmile="3\msx@60
\mathchardef\smallfrown="3\msx@61
\mathchardef\Subset="3\msx@62
\mathchardef\Supset="3\msx@63
\mathchardef\Cup="2\msx@64

\mathchardef\Cap="2\msx@65

\mathchardef\curlywedge="2\msx@66
\mathchardef\curlyvee="2\msx@67
\mathchardef\leftthreetimes="2\msx@68
\mathchardef\rightthreetimes="2\msx@69
\mathchardef\subseteqq="3\msx@6A
\mathchardef\supseteqq="3\msx@6B
\mathchardef\bumpeq="3\msx@6C
\mathchardef\Bumpeq="3\msx@6D
\mathchardef\lll="3\msx@6E

\mathchardef\ggg="3\msx@6F

\mathchardef\circledS="0\msx@73
\mathchardef\pitchfork="3\msx@74
\mathchardef\dotplus="2\msx@75
\mathchardef\backsim="3\msx@76
\mathchardef\backsimeq="3\msx@77
\mathchardef\complement="0\msx@7B
\mathchardef\intercal="2\msx@7C
\mathchardef\circledcirc="2\msx@7D
\mathchardef\circledast="2\msx@7E
\mathchardef\circleddash="2\msx@7F
\def\ulcorner{\delimiter"4\msx@70\msx@70 }
\def\urcorner{\delimiter"5\msx@71\msx@71 }
\def\llcorner{\delimiter"4\msx@78\msx@78 }
\def\lrcorner{\delimiter"5\msx@79\msx@79 }
\def\yen{\mathhexbox\msx@55 }
\def\checkmark{\mathhexbox\msx@58 }
\def\circledR{\mathhexbox\msx@72 }
\def\maltese{\mathhexbox\msx@7A }
\mathchardef\lvertneqq="3\msy@00
\mathchardef\gvertneqq="3\msy@01
\mathchardef\nleq="3\msy@02
\mathchardef\ngeq="3\msy@03
\mathchardef\nless="3\msy@04
\mathchardef\ngtr="3\msy@05
\mathchardef\nprec="3\msy@06
\mathchardef\nsucc="3\msy@07
\mathchardef\lneqq="3\msy@08
\mathchardef\gneqq="3\msy@09
\mathchardef\nleqslant="3\msy@0A
\mathchardef\ngeqslant="3\msy@0B
\mathchardef\lneq="3\msy@0C
\mathchardef\gneq="3\msy@0D
\mathchardef\npreceq="3\msy@0E
\mathchardef\nsucceq="3\msy@0F
\mathchardef\precnsim="3\msy@10
\mathchardef\succnsim="3\msy@11
\mathchardef\lnsim="3\msy@12
\mathchardef\gnsim="3\msy@13
\mathchardef\nleqq="3\msy@14
\mathchardef\ngeqq="3\msy@15
\mathchardef\precneqq="3\msy@16
\mathchardef\succneqq="3\msy@17
\mathchardef\precnapprox="3\msy@18
\mathchardef\succnapprox="3\msy@19
\mathchardef\lnapprox="3\msy@1A
\mathchardef\gnapprox="3\msy@1B
\mathchardef\nsim="3\msy@1C
\mathchardef\ncong="3\msy@1D

\mathchardef\varsubsetneq="3\msy@20
\mathchardef\varsupsetneq="3\msy@21
\mathchardef\nsubseteqq="3\msy@22
\mathchardef\nsupseteqq="3\msy@23
\mathchardef\subsetneqq="3\msy@24
\mathchardef\supsetneqq="3\msy@25
\mathchardef\varsubsetneqq="3\msy@26
\mathchardef\varsupsetneqq="3\msy@27
\mathchardef\subsetneq="3\msy@28
\mathchardef\supsetneq="3\msy@29
\mathchardef\nsubseteq="3\msy@2A
\mathchardef\nsupseteq="3\msy@2B
\mathchardef\nparallel="3\msy@2C
\mathchardef\nmid="3\msy@2D
\mathchardef\nshortmid="3\msy@2E
\mathchardef\nshortparallel="3\msy@2F
\mathchardef\nvdash="3\msy@30
\mathchardef\nVdash="3\msy@31
\mathchardef\nvDash="3\msy@32
\mathchardef\nVDash="3\msy@33
\mathchardef\ntrianglerighteq="3\msy@34
\mathchardef\ntrianglelefteq="3\msy@35
\mathchardef\ntriangleleft="3\msy@36
\mathchardef\ntriangleright="3\msy@37
\mathchardef\nleftarrow="3\msy@38
\mathchardef\nrightarrow="3\msy@39
\mathchardef\nLeftarrow="3\msy@3A
\mathchardef\nRightarrow="3\msy@3B
\mathchardef\nLeftrightarrow="3\msy@3C
\mathchardef\nleftrightarrow="3\msy@3D
\mathchardef\divideontimes="2\msy@3E
\mathchardef\varnothing="0\msy@3F
\mathchardef\nexists="0\msy@40
\mathchardef\mho="0\msy@66
\mathchardef\eth="0\msy@67
\mathchardef\eqsim="3\msy@68
\mathchardef\beth="0\msy@69
\mathchardef\gimel="0\msy@6A
\mathchardef\daleth="0\msy@6B
\mathchardef\lessdot="3\msy@6C
\mathchardef\gtrdot="3\msy@6D
\mathchardef\ltimes="2\msy@6E
\mathchardef\rtimes="2\msy@6F
\mathchardef\shortmid="3\msy@70
\mathchardef\shortparallel="3\msy@71
\mathchardef\smallsetminus="2\msy@72
\mathchardef\thicksim="3\msy@73
\mathchardef\thickapprox="3\msy@74
\mathchardef\approxeq="3\msy@75
\mathchardef\succapprox="3\msy@76
\mathchardef\precapprox="3\msy@77
\mathchardef\curvearrowleft="3\msy@78
\mathchardef\curvearrowright="3\msy@79
\mathchardef\digamma="0\msy@7A
\mathchardef\varkappa="0\msy@7B
\mathchardef\hslash="0\msy@7D
\mathchardef\hbar="0\msy@7E
\mathchardef\backepsilon="3\msy@7F
\def\Bbb{\relaxnext@\ifmmode\let\next\Bbb@\else
 \def\next{\Err@{Use \string\Bbb\space only in math mode}}\fi\next}
\def\Bbb@#1{{\Bbb@@{#1}}}
\def\Bbb@@#1{\noaccents@\fam\msyfam#1}

\catcode`\@=12

\font\tenfrak=eufm10
\font\sevenfrak=eufm7
\font\fivefrak=eufm5
\newfam\frakfam \def\frak{\fam\frakfam\tenfrak} \textfont\frakfam=\tenfrak
  \scriptfont\frakfam=\sevenfrak  \scriptscriptfont\frakfam=\fivefrak

\begin{document}
    \pagestyle{plain}
    \setlength{\baselineskip}{1.3\baselineskip}
    \setlength{\parindent}{\parindent}

\title{{\bf Lie algebroids associated to Poisson actions}}
\author{Jiang-Hua Lu\\
Department of Mathematics, University of Arizona, Tucson, AZ 85721}
\maketitle

\newtheorem{thm}{Theorem}[section]
\newtheorem{lem}[thm]{Lemma}
\newtheorem{prop}[thm]{Proposition}
\newtheorem{cor}[thm]{Corollary}
\newtheorem{rem}[thm]{Remark}
\newtheorem{exam}[thm]{Example}
\newtheorem{nota}[thm]{Notation}
\newtheorem{dfn}[thm]{Definition}
\newtheorem{ques}[thm]{Question}
\newtheorem{eq}{thm}

\newcommand{\rw}{\rightarrow}
\newcommand{\lrw}{\longrightarrow}
\newcommand{\rhu}{\rightharpoonup}
\newcommand{\lhu}{\leftharpoonup}
\newcommand{\Map}{\longmapsto}
\newcommand{\qed}{\begin{flushright} {\bf Q.E.D.}\ \ \ \ \
                  \end{flushright} }
\newcommand{\beqa}{\begin{eqnarray*}}
\newcommand{\eeqa}{\end{eqnarray*}}

\newcommand{\la}{\mbox{$\langle$}}
\newcommand{\ra}{\mbox{$\rangle$}}
\newcommand{\ot}{\mbox{$\otimes$}}
\newcommand{\xa}{\mbox{$x_{(1)}$}}
\newcommand{\xb}{\mbox{$x_{(2)}$}}
\newcommand{\xc}{\mbox{$x_{(3)}$}}
\newcommand{\ya}{\mbox{$y_{(1)}$}}
\newcommand{\yb}{\mbox{$y_{(2)}$}}
\newcommand{\yc}{\mbox{$y_{(3)}$}}
\newcommand{\yd}{\mbox{$y_{(4)}$}}
\renewcommand{\aa}{\mbox{$a_{(1)}$}}
\newcommand{\ab}{\mbox{$a_{(2)}$}}
\newcommand{\ac}{\mbox{$a_{(3)}$}}
\newcommand{\ad}{\mbox{$a_{(4)}$}}
\newcommand{\ba}{\mbox{$b_{(1)}$}}
\newcommand{\bt}{\mbox{$b_{(2)}$}}
\newcommand{\bc}{\mbox{$b_{(3)}$}}
\newcommand{\ca}{\mbox{$c_{(1)}$}}
\newcommand{\cb}{\mbox{$c_{(2)}$}}
\newcommand{\cc}{\mbox{$c_{(3)}$}}

\newcommand{\ts}{\mbox{$\sigma$}}
\newcommand{\las}{\mbox{${}_{\sigma}\!A$}}
\newcommand{\lasone}{\mbox{${}_{\sigma'}\!A$}}
\newcommand{\ras}{\mbox{$A_{\sigma}$}}
\newcommand{\rds}{\mbox{$\cdot_{\sigma}$}}
\newcommand{\lds}{\mbox{${}_{\sigma}\!\cdot$}}

\newcommand{\bb}{\mbox{$\bar{\beta}$}}
\newcommand{\bg}{\mbox{$\bar{\gamma}$}}

\newcommand{\id}{\mbox{${\em id}$}}
\newcommand{\Fun}{\mbox{${\em Fun}$}}
\newcommand{\End}{\mbox{${\em End}$}}
\newcommand{\Hom}{\mbox{${\em Hom}$}}
\newcommand{\ta}{\mbox{${\mbox{$\scriptscriptstyle A$}}$}}
\newcommand{\ms}{\mbox{${\mbox{$\scriptscriptstyle M$}}$}}
\newcommand{\ap}{\mbox{$A_{\mbox{$\scriptscriptstyle P$}}$}}
\newcommand{\tx}{\mbox{$\mbox{$\scriptscriptstyle X$}$}}
\newcommand{\pp}{\mbox{$\pi_{\mbox{$\scriptscriptstyle P$}}$}}
\newcommand{\pg}{\mbox{$\pi_{\mbox{$\scriptscriptstyle G$}}$}}
\newcommand{\asemi}{\mbox{$\ap \#_{\sigma} A^*$}}
\newcommand{\dsemi}{\mbox{$A \#_{\Delta} A^*$}}

\newcommand{\semi}{\mbox{$\times_{{\frac{1}{2}}}$}}
\newcommand{\fd}{\mbox{${\frak d}$}}
\newcommand{\fg}{\mbox{${\frak g}$}}
\newcommand{\fh}{\mbox{${\frak h}$}}
\newcommand{\fn}{\mbox{${\frak n}$}}
\newcommand{\fbp}{\mbox{${\frak b}_{+}$}}
\newcommand{\fbm}{\mbox{${\frak b}_{-}$}}
\newcommand{\fnp}{\mbox{${\frak n}_{+}$}}
\newcommand{\fnm}{\mbox{${\frak n}_{-}$}}
\newcommand{\fgs}{\mbox{${\frak g}^*$}}
\newcommand{\wg}{\mbox{$\wedge {\frak g}$}}
\newcommand{\wgs}{\mbox{$\wedge {\frak g}^*$}}
\newcommand{\wxl}{\mbox{$x_1 \wedge x_2 \wedge \cdots \wedge x_l$}}
\newcommand{\wxk}{\mbox{$x_1 \wedge x_2 \wedge \cdots \wedge x_k$}}
\newcommand{\wyl}{\mbox{$y_1 \wedge y_2 \wedge \cdots \wedge y_l$}}
\newcommand{\wxkm}{\mbox{$x_1 \wedge x_2 \wedge \cdots \wedge x_{k-1}$}}
\newcommand{\wxik}{\mbox{$\xi_1 \wedge \xi_2 \wedge \cdots \wedge \xi_k$}}
\newcommand{\wxikm}{\mbox{$\xi_1 \wedge \cdots \wedge \xi_{k-1}$}}
\newcommand{\wetal}{\mbox{$\eta_1 \wedge \eta_2 \wedge \cdots \wedge \eta_l$}}

\newcommand{\winv}{\mbox{$(\wedge \fg_{1}^{\perp})^{\fg_1}$}}
\newcommand{\wetak}{\mbox{$\eta_1 \wedge \cdots \wedge \eta_k$}}
\newcommand{\gonep}{\mbox{$\fg_{1}^{\perp}$}}
\newcommand{\wonep}{\mbox{$\wedge \fg_{1}^{\perp}$}}

\newcommand{\db}{\mbox{$\fd = \fg \bowtie \fgs$}}
\newcommand{\fds}{\mbox{${\scriptscriptstyle {\frak d}}$}}
\newcommand{\fl}{\mbox{${\frak l}$}}

\newcommand{\Gs}{\mbox{$G^*$}}
\newcommand{\pis}{\mbox{$\pi_{\sigma}$}}
\newcommand{\ea}{\mbox{$E_{\alpha}$}}
\newcommand{\eb}{\mbox{$E_{-\alpha}$}}
\newcommand{\Bm}{\mbox{$ {}^B \! M$}}
\newcommand{\kBm}{\mbox{$ {}^B \! M^k$}}
\newcommand{\Bb}{\mbox{$ {}^B \! b$}}
\renewcommand{\epsilon}{\mbox{$\varepsilon$}}

\newcommand{\cfg}{\mbox{$C(\fg \oplus \fgs)$}}
\newcommand{\ps}{\mbox{$\pi^{\#}$}}
\newcommand{\backl}{\mathbin{\vrule width1.5ex height.4pt\vrule height1.5ex}}

\newcommand{\bx}{\mbox{${\bar{x}}$}}
\newcommand{\by}{\mbox{${\bar{y}}$}}
\newcommand{\bz}{\mbox{${\bar{z}}$}}
\newcommand{\pgs}{\mbox{${\pi_{\mbox{\tiny G}^{*}}}$}}

\newcommand{\tp}{\mbox{$\varphi$}}
\newcommand{\sn}{\mbox{$s_{\scriptscriptstyle N}$}}
\newcommand{\tn}{\mbox{$t_{\scriptscriptstyle N}$}}
\newcommand{\sm}{\mbox{$s_{\scriptscriptstyle M}$}}
\newcommand{\tm}{\mbox{$t_{\scriptscriptstyle M}$}}
\newcommand{\en}{\mbox{$\epsilon_{\scriptscriptstyle N}$}}
\newcommand{\mem}{\mbox{$\epsilon_{\scriptscriptstyle M}$}}

\section{Introduction}
\label{sec_intro}

This work is motivated by a result of Drinfeld in \cite{dr:homog}. Recall
\cite{dr:bigbra}
that a {\bf Poisson Lie group} is a Lie group $G$ together with a Poisson
structure such that the group multiplication map
\[
G \times G \lrw G
\]
is a Poisson map. Given a Poisson Lie group $G$ and a Poisson manifold
$P$, an action
\[
\sigma: ~ G \times P \lrw P
\]
of $G$ on $P$ is called a
{\bf Poisson action} if the action map $\sigma$ is a Poisson
map. When the action is transitive, we say that
$P$ is a {\bf Poisson homogeneous $G$-space}.
Poisson $G$-spaces are the semi-classical
analogs of quantum spaces with
quantum group actions.
Special cases of Poisson homogeneous $G$-spaces can be
found in \cite{daso:affine}
\cite{lu:thesis} \cite{za:homog}.

\bigskip
Let $P$ be a Poisson homogeneous $G$-space. In \cite{dr:homog}, Drinfeld
shows that corresponding to each $p \in P$, there is a maximal isotropic Lie
subalgebra ${\frak l}_p$ of the Lie algebra $\fd$, the double Lie algebra
of the tangent Lie bialgebra $(\fg, \fgs)$ of $G$. Moreover, for $g \in G$, the
two Lie algebras ${\frak l}_p$ and ${\frak l}_{gp}$ are related by
${\frak l}_{gp} = Ad_{g}{\frak l}_p$ via the Adjoint action of $G$ on $\fd$.
In particular, they are isomorphic as Lie algebras. The Lie algebra
${\frak l}_p$ determines the Poisson structure on $P$,
and it can be used to classify Poisson homogeneous $G$-spaces \cite{dr:homog}.

\bigskip
The purpose of this note is to find an invariant setting for these Lie algebra
${\frak l}_p$'s. We construct, for every Poisson manifold $P$ with a
Poisson $G$-action (not necessarily $G$-homogeneous), a
Lie algebroid structure on
the direct sum vector bundle $(P \times \fg) \oplus T^*P$ over $P$.
This Lie
algebroid will be denoted by $A = (P \times \fg) \bowtie
T^*P$ to indicate the fact that it is built out of the two
Lie algebroids $P \times \fg$ and $T^*P$ and a pair of representations
of them on each other. Moreover, the Lie algebroid $A$ is
naturally equipped with an action of $G$, making it into a
Harish-Chandra Lie algebroid
of $G$ \cite{bb:jantzen}. It follows
immediately that the kernel
${\frak l}_p$ of the anchor map of $A$
at each $p \in P$
has a natural Lie algebra structure, and ${\frak l}_{gp} = g \cdot {\frak
l}_p$,
where $g \cdot {\frak l}_p$ comes from the action of $G$ on $A$.

\bigskip
When the anchor map of $A$ has full rank everywhere (it is said to
be transitive in such a case),
the subbundle of
$A$ defined by the kernel of the anchor map is a Lie
algebra bundle in the sense that local trivilizations exist \cite{mk:book}.
This is the case when $P$ is $G$-homogeneous or when $P$ is symplectic.
In both cases,
each Lie algebra ${\frak l}_p$ can be naturally
embedded in the Lie algebra $\fd$, the double of the tangent Lie bialgebra
of $G$, as a maximal isotropic Lie subalgebra.

When $P$ is $G$-homogeneous, these maximal isotropic Lie subalgebras
of $\fd$ are precisely
those described in \cite{dr:homog}.

When $P$ is symplectic, we show that the Lie algebra ${\frak l}_p$'s
define an isomorphic family of Lie bialgebras over $\fgs$.

\bigskip
As further applications, we describe the symplectic leaves of a
Poisson homogeneous $G$-space $P$ in terms of the Lie
algebra ${\frak l}_p$.  We show that the $G$-invariant
Poisson cohomology of $P$ can be realized as relative Lie
algebra cohomology of ${\frak l}_p$ relative to the
stabilizer subgroup of $G$ at $p$. As
an example, we calculate the $K$-invariant Poisson
cohomology of the Bruhat Poisson structure \cite{lu-we:poi}
on the generalized flag manifold $K/T$, where $K$ is a
compact semisimple Lie group and $T \subset K$ is a maximal
torus of $K$, and we show that there is exactly one cohomology class
for each element in the Weyl group of $K$.
We also describe the Lie groupoid
corresponding to the Lie algebroid $A = (P \times \fg) \bowtie
T^*P$ and show that it is an example of a double Lie groupoid
in the sense of Mackenzie \cite{mk:double}.

\bigskip
{\bf Acknowledgement} The author would like to thank Yvette Kosmann
Schwarzbach for drawing her attention to \cite{dr:homog} and Sam
Evens for answering many questions. She also would like to thank
Kirill Mackenzie for helpful comments.

\section{Lie algebroids}
\label{sec_lie-algebroids}

In this section, we collect some relevant facts about Lie algebroids.  More
systematic treatments of them can be found in \cite{mk:book} \cite{h-mk:cat}
and \cite{xu-mk:bi}.

\bigskip
A {\bf Lie algebroid} over a smooth manifold $P$ is a vector bundle
$A$ over $P$ together with a Lie algebra structure on
the space $\Gamma(A)$ of smooth sections of $A$ and  a
bundle map $\rho: A \rightarrow TP$ (called the anchor map of the
Lie algebroid) such that

1) $\rho$
defines a Lie algebra homomorphism  from $\Gamma(A)$ to the
space $\chi(P)$ of vector fields with the commutator
Lie algebra structure, and

2) for $f \in C^{\infty}(P), \omega_{1},
\omega_{2} \in \Gamma(A)$ the  following derivation law holds:
\[
\{\omega_1, f \omega_2 \} = f \{\omega_1, \omega_2\} + (\rho (\omega_1)
f) \omega_2.
\]
Immediate examples of Lie algebroids are 1) the tangent bundle $TP$
as a Lie algebroid over $P$ with the identity map as the anchor map, and 2)
a Lie algebra $\fg$ as a Lie algebroid over a one point space.

\bigskip
\noindent
{\bf Transformation Lie algebroids}.
Let $\sigma: ~ G \times P \lrw P $
be an action of a Lie group $G$ on a manifold $P$.
For each $x \in \fg$, the Lie algebra of $G$,
denote by $\sigma_x$ the vector field on $P$ defined by
\[
\sigma_x (p) ~ = ~ {\frac{d}{dt}}|_{t=0} \sigma(\exp tx, p), ~~~ p \in P.
\]
Then there is a natural Lie algebroid structure on the
trivial vector bundle $P \times \fg$ with
$-\sigma$ as the anchor map. Here we regard $\sigma: x \mapsto \sigma_x$
as a bundle map from $P \times \fg$ to $TP$.
The Lie bracket on the
space $\Gamma(P \times \fg) \cong C^{\infty}(P, \fg)$
of smooth sections of $P \times \fg$ is
given by
\begin{equation}
\label{eq_transformation-algebroid}
\{ \bx, ~ \by\} ~ = ~ [\bx, ~ \by]_{\fg} ~ - ~ \sigma_{\bar{x}} \cdot \by ~ + ~
\sigma_{\bar{y}} \cdot \bx,
\end{equation}
where
the first term on the right hand side denotes the pointwise Lie bracket
in ${\frak g}$, and the second term denotes the
derivative of the ${\frak g}$-valued function $\by$ in the direction of the
vector field $\sigma_{\bar{x}}$.

\bigskip
\noindent
{\bf Cotangent bundle Lie algebroids.}
Let $(P, \pi)$ be a Poisson manifold with Poisson bi-vector field $\pi$.
We use $\ps$ to denote the bundle map
\[
\ps(p): ~~T_{p}^{*}P \longrightarrow T_pP:~~
\alpha_p \Map - \alpha_p \backl \pi(p),
\]
or
\[
(\beta_p, ~ \ps(\alpha_p)) ~ = ~ \pi(p) (\beta_p, ~ \alpha_p), ~~~~~ \alpha_p,
\beta_p \in T_{p}^{*}P.
\]
(Note the sign convention here). Then,
with $-\ps$ as the anchor map, the cotangent bundle $T^*P$
becomes a Lie algebroid over $P$, where the Lie bracket on the space
$ {\Omega}^{1} ( P )$ of  1-forms on
$P$ is given by
\begin{eqnarray}
\label{eq_bracket-on-one-forms}
\{\alpha, \beta\} & = &d \pi (\alpha , \beta) - \pi^{\#}
\alpha \backl d \beta + \pi^{\#} \beta \backl d \alpha \\
\label{eq_2}
& = & - d \pi (\alpha , \beta) - L_{\pi^{\#} \alpha} \beta +
L_{\pi^{\#} \beta} \alpha,
\end{eqnarray}
where $L_{\pi^{\#} \alpha} \beta$ denotes the Lie derivative of
the $1$-form $\beta$ in the direction of the vector field $\pi^{\#} \alpha$.

\bigskip
\noindent
{\bf Lie algebroid morphisms.} The definition of Lie algebroid morphisms
between Lie algebroids over different bases is rather involved
\cite{h-mk:cat} \cite{xu-mk:bi}, the reason being that a bundle map
does not necessarily induce a map between sections. We will only need
the definition for the following special case.
Let $A$ be a Lie algebroid over $P$ with anchor map
$\rho$. Let $\fg$ be a Lie algebra, considered as a Lie
algebroid over a one point space. A smooth map
\[
\phi: ~ A \lrw \fg
\]
which is linear on each fiber is said to be a Lie algebroid morphism
if for any sections $\omega_1$ and $\omega_2$ of $A$,
\[
\phi \{\omega_1, ~ \omega_2\} ~ = ~ [\phi(\omega_1), ~ \phi(\omega_2)]_{{\frak
g}}
{}~ + ~ \rho(\omega_1) \cdot \phi(\omega_2)
{}~ - ~ \rho(\omega_2) \cdot \phi(\omega_1),
\]
where the first term on the right hand side denotes the pointwise
Lie bracket in $\fg$, and the second term denotes the Lie derivative
of the $\fg$-valued function $\phi(\omega_2)$ in the
direction of the vector field $\rho(\omega_1)$.

\begin{exam}
\label{exam_to-fg}
{\em
In the case of a transformation Lie algebroid $A = P \times \fg$, the
projection map $P \times \fg \rightarrow \fg$ is a Lie algebroid morphism.
}
\end{exam}

\bigskip
\noindent
{\bf Representations of Lie algebroids.} Let $A$ be a Lie algebroid
over $P$ with anchor map $\rho$. Let $E$ be a vector bundle over $P$.
A representation of $A$ on $E$ is a $k$-linear map
\[
\Gamma(A) \ot \Gamma(E) \lrw \Gamma(E): ~~ a \ot s \Map D_a s,
\]
where $\Gamma(A)$ and $\Gamma(E)$ denote respectively the
spaces of smooth sections of $A$ and $E$, such that for
any $a, b \in \Gamma(A), s \in \Gamma(E)$ and $f \in C^{\infty}(P)$,
\beqa
& & (1) ~~ D_{fa}(s) ~ = ~ f D_a s;\\
& & (2) ~~ D_a(fs) ~ = ~ f D_a s + (\rho(a)f) s;\\
& &  (3) ~~ D_a(D_bs) - D_b(D_as) ~ = ~ D_{\{a, b\}} s.
\eeqa
If $\la ~,~\ra$ is a smooth section of $S^2(E^*)$, the second
symmetric power of the dual bundle of $E$, we say that the representation
of $A$ on $E$ is $\la ~,~\ra$-preserving if for any $a \in
\Gamma(A), ~ s_1, s_2 \in \Gamma(E)$,
\[
\rho(a) (s_1, ~ s_2) ~ = ~ \la D_a s_1, ~ s_2 \ra ~ + ~
\la s_1, ~ D_a s_2 \ra.
\]

\begin{exam}
\label{exam_rep1}
{\em
A representation of the tangent bundle $TP$ over $P$ is a vector bundle
$E$ over $P$ together with a flat linear connection.
}
\end{exam}

\begin{exam}
\label{exam_rep2}
{\em
Any transformation Lie algebroid $P \times \fg$ has a natural
representation on the tangent bundle
$TP$ of $P$ via
\begin{equation}
\label{eq_trans-on-TP}
D_{\bar{x}}V ~ = ~ -[\sigma_{\bar{x}}, ~ V] ~ - ~ \sigma_{V \cdot {\bar{x}}},
\end{equation}
where $V$ is a vector field on $P, ~ \bx$ is a section of
$P \times \fg$, considered as a $\fg$-valued function of $P$, and
$V \cdot \bx$ denotes the Lie derivative of $\bx$ in the direction of $V$.
Correspondingly, there a representation of $P \times \fg$ on
the cotangent bundle $T^*P$ satisfying
\begin{equation}
\label{eq_trans-on-T*P-1}
- \sigma_{\bar{x}} (\alpha, ~ V) ~ = ~ (D_{\bar{x}} \alpha, ~ V)
{}~ + ~ (\alpha, ~ D_{\bar{x}} V)
\end{equation}
for any $1$-form $\alpha$ on $P$.
Or, equivalently,
\begin{equation}
\label{eq_trans-on-T*P-2}
(D_{\bar{x}} \alpha, ~ V) ~ = ~ - \sigma_{\bar{x}} (\alpha, ~ V)
{}~ + ~ (\alpha, ~ [\sigma_{\bar{x}}, ~ V] ~ + ~ \sigma_{V \cdot {\bar{x}}}).
\end{equation}
Notice that when $\bar{x}$ is a constant section of $P \times \fg$
corresponding to $x \in \fg$, the $1$-form $D_{x} \alpha$ is given by
the Lie derivative
\[
D_x \alpha ~ = ~ - L_{\sigma_x} \alpha.
\]
}
\end{exam}

\begin{exam}
\label{exam_trans-on-trivial}
{\em
Let $A$ be a Lie algebroid over $P$ with anchor map $\rho$. Suppose
that $\phi: A \rightarrow \fg$
is a Lie algebroid morphism from $A$ to a Lie algebra $\fg$. Then for
any vector space $U$ with a $\fg$-action, there is a
representation of $A$ on the trivial vector bundle
$P \times U$ given by
\begin{equation}
\label{eq_trans-on-U}
D_a \bar{u} ~ = ~ \rho(a) \cdot \bar{u} ~ + ~ \phi(a) (\bar{u}),
\end{equation}
where $a$ is a section of $A, ~ \bar{u}$ is a section of $P \times U$,
and $\phi(a) (\bar{u})$ denotes the action of $\phi(a) \in C^{\infty}(P,
\fg)$ on $\bar{u} \in C^{\infty}(P, U)$ taken pointwise over $P$. In
particular, there
is a representation of $A$ on the
trivial vector bundle $P \times \fgs$ given by
\[
D_a \bar{\xi} ~ = ~ \rho(a) \cdot \bar{\xi} ~ + ~ ad_{\phi(a)}^{*}
\bar{\xi},
\]
where $\bar{\xi} \in C^{\infty}(P, \fgs)$ is a
section of $P \times \fgs$, and the co-adjoint representation of
$\fg$ on $\fgs$ is defined by
\[
(ad_{x}^{*} \xi, ~ y) ~ = ~ -(\xi, ~ [x, y]), ~~~~~~x, y \in \fg, ~ \xi \in
\fgs.
\]
}
\end{exam}

\bigskip
\noindent
{\bf Lie algebroid cohomology} Let $A$ be a Lie algebroid over $P$
with anchor map $\rho$. Let $E$ be a representation of $A$.
Define
\[
d: ~ \Gamma(\wedge^{k-1} A^* \otimes E) \lrw \Gamma(\wedge^{k} A^*
\otimes E), ~~~k =1, 2, 3, ...
\]
by
\begin{equation}
\label{eq_dk}
df (a_1, ..., a_k)  =  \sum_i (-1)^{i+1} D_{a_i} f(a_1, ..., \hat{a}_i,
..., a_k)  +  \sum_{i < j} (-1)^{i+j} f(\{a_i, a_j\}, ..., \hat{a}_i, ...,
\hat{a}_j, ..., a_k).
\end{equation}
Then $d$ is well-defined and that $d^2 = 0$. The
cohomology of $(C^{\bullet}(A, E), d)$ is called the Lie algebroid cohomology
of $A$ with coefficients in $E$, and it is denoted by $H^{\bullet}(A, E)$.
Elements of $H^{0}(A,E)$ are called $A$-parallel sections of $E$.
When $E$ is the trivial $1$-dimensional
vector bundle with the representation of $A$ given by
\[
D_a f ~ = ~ \rho(a) f, ~~~~ f \in C^{\infty}(P),
\]
the cohomology of $A$ with coefficients in $E$ is
called the cohomology of $A$ with trivial coefficients.

In the case when $A = TP$ is the tangent bundle of $P$, the Lie algebroid
cohomology of $A$ with trivial coefficients is nothing but
the de Rham cohomology of $P$.

In the case when $A = T^*P$ is the cotangent bundle Lie algebroid
of a Poisson manifold $(P, \pi)$, the Lie algebroid cohomology of
$A$ with trivial coefficients is called the Poisson cohomology
of $(P, \pi)$. The cochain complex
is the space $\Gamma(\wedge^{\bullet} TP)$ of multi-vector
fields on $P$, and the coboundary operator $d$ can be expressed as
$d = [\pi, ~ \bullet]$, where $[ ~ , ~ ]$
stands for the
Schouten bracket \cite{ku:schouten} on the multi-vector fields on $P$
\cite{li:poi}

In the case when $A = P \times \fg$ is a transformation Lie algebroid, the
Lie algebra cohomology of $A$ with trivial coefficients is simply
the Lie algebra cohomology of $\fg$ with coefficients in $C^{\infty}(P)$.

\bigskip
\noindent
{\bf Transitive Lie algebroids}. A Lie algebroid $A$ over a manifold $P$
with anchor map $\rho: A \rightarrow TP$ is said to be transitive if
$\rho$ has full rank at every point of $P$. In this case, the subbundle
$L$ of $A$ defined by the kernel of $\rho$ is a Lie algebra bundle
in the sense that each fiber has a Lie algebra structure and
local trivilizations exist \cite{mk:book}. Any section
$\gamma: TP \rightarrow A$ of $\rho$
defines a linear connection $\nabla$
of $L$ by
\[
\nabla_{V}^{\gamma} (l) ~ = ~ \{ \gamma (V), ~~ l\}, ~~ ~ V \in \chi(P), ~~
l \in \Gamma(L),
\]
and it satisfies
\begin{equation}
\label{eq_L}
\nabla_{V}^{\gamma} [l_1, ~ l_2] ~ = ~ [\nabla_{V}^{\gamma} l_1, ~ l_2]
{}~ + ~ [l_1, ~ \nabla_{V}^{\gamma} l_2], ~~~~\forall l_1, ~ l_2 \in
\Gamma(L).
\end{equation}

\section{Poisson Lie groups and Poisson actions}
\label{sec_poi}

In this section, we
collect some facts about Poisson Lie groups that will
be used in this paper. Details can be found in
\cite{dr:bigbra} \cite{sts:dressing} \cite{ks:podr} and \cite{lu:thesis}.

\bigskip
Assume that $(G, \pg)$ is a Poisson Lie group with
Poisson bi-vector field $\pg$.
For each $g \in G$,
the element $r_{g^{-1}} \pg(g)$ is in $\wedge^2 \fg$, where $r_{g^{-1}}$
denotes the right translation on $\wedge^2 TG$ by $g^{-1}$.
The derivative of the map
\[
G \lrw \wedge^2 \fg: ~~ g \Map r_{g^{-1}} \pg(g)
\]
at the identity element $e$ of $G$ is denoted by $\delta$, so we have
\begin{equation}
\label{eq_delta}
\delta: ~~ \fg \lrw \wedge^2 \fg: ~~ \delta(x) ~ = ~ {\frac{d}{dt}}|_{t=0}
r_{\exp(-tx)} \pg (\exp tx).
\end{equation}
The dual map of $\delta$ defines a Lie bracket
on $\fgs$, i.e.,
\begin{equation}
\label{eq_bra-fgs}
([\xi, ~ \eta], ~ x) ~ = ~ (\delta(x), ~~ \xi \wedge \eta), ~~~~~~~
x \in \fg, ~ \xi, \eta \in \fgs.
\end{equation}
The pair $(\fg, \delta)$ (or the pair $(\fg, \fgs)$)
is called the tangent  Lie
bialgebra of the Poisson Lie group $G$ \cite{dr:bigbra}.
For $x \in \fg$ and $\xi \in \fgs$, denote by
$ad_{x}^{*} \xi \in \fgs$ and $ad_{\xi}^{*} x \in \fg$
respectively the elements given by
\beqa
& & \la ad_x^* \xi, ~ y \ra ~ = ~ -\la \xi, ~ [x, y] \ra  ~~~~~
\forall y \in \fg\\
& & \la ad_{\xi}^* x, ~ \eta \ra ~ = ~ - \la x, ~ [\xi, \eta] \ra ~~~~~
\forall \eta \in \fgs.
\eeqa
The double Lie algebra $\fd$ is the vector space
$\fg \oplus \fgs$ with the Lie bracket
\begin{equation}
\label{eq_bra-on-d}
[x_1 + \xi_1, ~ x_2 + \xi_2] ~ = ~ [x_1, ~ x_2] ~ + ~ ad_{\xi_1}^{*} x_2 ~
- ~ ad_{\xi_2}^{*} x_1 ~ + ~ [\xi_1, ~ \xi_2] ~ + ~
ad_{x_1}^* \xi_2 ~ - ~ ad_{x_2}^{*} \xi_1.
\end{equation}
It is the unique Lie bracket on the vector space $\fg \oplus \fgs$
characterized
by the properties that 1)
when restricted to $\fg$
and $\fgs$, it coincides with the given Lie brackets on $\fg$ and $\fgs$,
and 2) the scalar product $\la ~~, ~~\ra_{\fd}$ on $\fg \oplus \fgs$
defined by
\begin{equation}
\label{eq_scalar}
\la x_1 + \xi_1, ~ x_2, + \xi_2 \ra_{\fd} ~ = ~ (x_1, ~ \xi_2) ~ + ~ (\xi_1 , ~
x_2), ~~~ x_1, ~ x_2 \in \fg, ~ \xi_1, ~ \xi_2 \in \fgs
\end{equation}
is ad-invariant with respect to
this Lie bracket. We use $\fd = \fg \bowtie \fgs$ to denote this Lie
algebra, indicating
the fact that the Lie bracket is built out of the Lie
brackets on $\fg$ and on $\fgs$ and their co-adjoint actions on each
other.

Let $G^*$ be the simply-connected Lie group with Lie algebra
$\fgs$. Then it is also a Poisson Lie group with tangent
Lie bialgebra $(\fgs, \fg)$. It is called the dual Poisson Lie
group of $G$.

Assume that $D$ is a connected Lie group with Lie algebra $\fd$.
Assume also that both $G$ and $G^*$ can be embedded in $D$
as subgroups corresponding to
the Lie algebra inclusions
$\fg \subset \fd$ and $\fgs \subset \fd$.
Assume furthermore that every element of $D$ can be uniquely written as
a product $gu$ for some $g \in G$ and $u \in G^*$. Then we
can identify $G$ with the quotient space $D/G^*$.
Consequently, there is a left action of $D$ on $G$
by left translations. We will denote it by
\begin{equation}
\label{eq_D-on-G-left}
D \times G \lrw G: ~~ (d, ~ g) \Map ^{d}\!\!g.
\end{equation}
Infinitesimally, this defines a Lie algebra anti-homomorphism
from $\fd$ to the Lie algebra $\chi(G)$ of vector fields
on $G$:
\begin{equation}
\label{eq_fd-on-G-left}
\fd \lrw \chi(G): ~ x + \xi \Map \rho_{x + \xi }(g) ~ = ~
r_g x ~ - ~
r_g \left( \xi \backl (r_{g^{-1}} \pg (g)) \right), ~~~~
g \in G.
\end{equation}
This action of $\fd$ on $G$ will be used in
Example \ref{exam_G}.
The vector fields $-\rho_{\xi}, ~ \xi \in \fgs$, are called the right
dressing vector fields on $G$ \cite{lu:thesis}.
When the global decomposition $D = GG^*$ does not hold,
we will assume that these vector fields integrate
to a left action of $D$ on $G$.

Similarly, by identifying $G$ with the quotient space
$G^* \backslash D$, we get a right action of $D$ on $G$ which we
denote by
\begin{equation}
\label{eq_D-on-G-right}
G \times D \lrw G: ~~ (g, ~ d) \Map g^d.
\end{equation}
The corresponding infinitesimal right action of $D$ on
$G$ is given by
\begin{equation}
\label{eq_fd-on-G-right}
\fd \lrw \chi(G): ~ x + \xi \Map \lambda_{x + \xi}(g) ~ = ~
l_g x ~ + ~ r_g(Ad_{g^{-1}}^{*} \xi \backl (r_{g^{-1}} \pg (g))).
\end{equation}
This action will be used in Theorem \ref{thm_leaves}.
The vector fields $-\lambda_{\xi}$,  for $\xi \in \fgs$, are called
the left dressing vector fields on $G$ \cite{lu:thesis}.

On the other hand, the group $G$ acts on $\fd$ by Lie algebra automorphisms via
\begin{equation}
\label{eq_G-on-fd}
Ad_g  (x + \xi) ~ = ~ Ad_g x + (Ad_{g^{-1}}^{*} \xi) \backl (r_{g^{-1}}
\pg(g)) ~ + ~ Ad_{g^{-1}}^{*} \xi.
\end{equation}
The corresponding action of $\fg$ on $\fd$ is simply the adjoint action
of $\fd$ on itself restricted to $\fg$. This action will be used
in Proposition \ref{prop_to-fd}. Similarly, the group
$G^*$ acts on $\fd$ by Lie algebra automorphisms via
\begin{equation}
\label{eq_Gs-on-fd}
Ad_u  (x + \xi) ~ = ~ Ad_{u^{-1}}^{*} x ~ + ~ Ad_{u} \xi ~ + ~
(Ad_{u^{-1}}^{*} x) \backl (r_{u^{-1}} \pi_{\mbox{\tiny G}^{*}}(u)).
\end{equation}
This action will be used in Examples \ref{exam_symplectic-iso}
and \ref{exam_dressing-orbit}.

\bigskip
We now turn to Poisson actions by $G$.
For any action
$\sigma: G \times
P \rightarrow P$
of $G$ on a manifold $P$, we can define
\begin{equation}
\label{eq_phi}
\phi: ~ T^*P \lrw \fgs: ~~ (\phi(\alpha_p), ~ x) ~ = ~ (\alpha_p, ~
\sigma_x(p))
,
{}~~~~ x \in \fg, ~ \alpha_p \in T_{p}^{*}P, ~ p \in P.
\end{equation}
Each $1$-form $\alpha$ on $P$ then defines a $\fgs$-valued function
$\phi(\alpha)$ on $P$:
\[
\phi(\alpha)(p) ~ = ~ \phi(\alpha_p).
\]

Assume now that $(P, \pi)$ is a Poisson manifold and $G$ is a
Poisson Lie group.
Recall that an action $\sigma: G \times
P \rightarrow P$ is said to be Poisson if
for any $g \in G, p \in P$
\[
\pi(gp) ~ = ~ g_* \pi(p) ~ + ~ p_* \pg(g),
\]
where $p_*$ is the differential of the map from $G$ to $P$
given by $g \mapsto gp$. When $G$ is connected, this is equivalent to the
following infinitesimal criterion \cite{lu-we:poi}:
\begin{equation}
\label{eq_poisson-action-infinitesimal}
L_{\sigma_{\bar{x}}} \pi ~ = ~ \sigma_{\delta(x)}
\end{equation}
for all $x \in \fg$, where $\delta(x) \in \fg \wedge \fg$ is given by
(\ref{eq_delta}), and $\sigma_{\delta(x)} = \sum_i \sigma_{x_i}
\wedge \sigma_{y_i}$ if $\delta(x) = \sum_i x_i \wedge y_i$.
This is then equivalent to
\begin{equation}
\label{eq_phi-liebi}
\phi(\{\alpha, \beta\}) ~ = ~ [\phi(\alpha), ~ \phi(\beta)]_{\fgs} ~
- ~ \ps(\alpha) \cdot \phi(\beta) ~ + ~ \ps(\beta) \cdot \phi(\alpha)
\end{equation}
for any $1$-forms $\alpha$ and $\beta$ on $P$. In other words,
we have the following Proposition:

\begin{prop} \cite{xu:poisson-groupoid}
\label{prop_phi-liebi}
Assume that $G$ is connected.
The action $\sigma$ is a Poisson action if and only if
the map
\[
\phi: ~ T^*P \lrw \fgs
\]
defined by
(\ref{eq_phi}) is a Lie algebroid morphism,
where $T^*P$ has the cotangent bundle Lie
algebroid structure defined by the Poisson structure $\pi$ on $P$, and
$\fgs$ has the Lie algebra structure given by
(\ref{eq_bra-fgs}), considered as a Lie algebroid
over a one point space.
\end{prop}

\section{The Lie algebroid $A = (P \times \fg) \bowtie T^*P$}
\label{sec_A}

In this section, we describe, for every Poisson manifold $P$ with a
Poisson action $\sigma: G \times P \rightarrow P$
of a Poisson Lie group $G$, a Lie algebroid structure over
$P$ on the direct sum vector bundle $(P \times \fg) \oplus T^*P$.

\bigskip
First, the action of $G$ on $P$ defines a representation of
the transformation Lie algebroid $P \times \fg$ on $T^*P$
given by
(\ref{eq_trans-on-T*P-2}) as
in Example \ref{exam_rep2}.

Secondly, using the Lie groupoid homomorphism
\[
\phi: ~ T^*P \lrw \fgs
\]
in Proposition \ref{prop_phi-liebi}, we get
(see Example \ref{exam_trans-on-trivial})
a representation of the cotangent bundle Lie
algebroid $T^*P$ on the trivial vector bundle $P \times \fg$
induced from the co-adjoint representation of $\fgs$ on $\fg$. It is given by
\begin{equation}
\label{eq_T*P-on-trans}
D_{\alpha} \bx ~ = ~ - \pi^{\#} (\alpha) \cdot \bx ~ + ~ ad_{\phi(\alpha)}^{*}
\bx
\end{equation}
where $(ad_{\phi(\alpha)}^{*} \bx)(p) = ad_{\phi(\alpha)(p)}^{*} \bx(p)$ for
$p \in P$.

\begin{thm}
\label{thm_A}
Let $\sigma: G \times P \rightarrow P$ be a Poisson action of the Poisson
Lie group $G$ on  the Poisson manifold $P$.
Let
\begin{equation}
\label{eq_dfn-A}
A ~ \stackrel{def}{=} ~ (P \times \fg) \oplus T^*P
\end{equation}
be the direct sum vector bundle over $P$. Then
there is a Lie algebroid structure on $A$ such that
both the transformation Lie algebroid $(P \times \fg, -\sigma)$ and
the cotangent bundle Lie algebroid $(T^*P, -\ps)$
are Lie subalgebroids of $A$. The anchor map of $A$ is $-\sigma - \ps$.
The Lie bracket between a section $\bx$ of $P \times \fg$ and a section
$\alpha$ of $T^*P$ is given by
\begin{equation}
\label{eq_bracket-on-A}
\{\bx, ~ \alpha\} ~ = ~  - D_{\alpha} \bx ~ + ~ D_{\bar{x}} \alpha,
\end{equation}
where $D_{\alpha} \bx$ is given by (\ref{eq_T*P-on-trans}) and
$D_{\bar{x}} \alpha$
by (\ref{eq_trans-on-T*P-2});
\end{thm}

\begin{dfn}
\label{dfn_A}
{\em
We use $A = (P \times \fg) \bowtie T^*P$ to denote this
Lie algebroid, indicating the fact that the Lie bracket on
the sections of $A$ uses both the representation of $P \times \fg$
on $T^*P$ and the representation of $T^*P$ on
$P \times \fg$.
}
\end{dfn}

Examples will be given in Section \ref{sec_exam}.
The proof of Theorem \ref{thm_A} will be given in Section \ref{sec_proof}.
We now discuss some properties of the Lie algebroid $A$.

\bigskip
By considering each element of $\fg$ as a constant section of
$A$, we get a Lie algebra homomorphism
\begin{equation}
\label{eq_i}
i: ~ \fg \lrw \Gamma(A): ~~ i(x)_p = x, ~~~~ \forall p \in P.
\end{equation}
It induces an action of $\fg$ on $\Gamma(A)$ by
\begin{equation}
\label{eq_fg-on-A}
\fg \ni x: ~ (\by, ~ \beta) \Map \{i(x), ~ (\by, ~ \beta)\} ~ = ~ \left(
[x, \by]_{\fg} - \sigma_{x} \cdot \by - ad_{\phi(\beta)}^{*} x, ~
{}~-L_{\sigma_x} \beta \right),
\end{equation}
where $\by \in C^{\infty}(P, \fg)$ and $\beta$ is a $1$-form on $P$.

\begin{prop}
\label{prop_harish-chandra-A}
The following defines a left action of $G$ on $A$:
\begin{equation}
\label{eq_G-on-A}
g \cdot(x, ~\alpha_p) ~ = ~ (Ad_g x + Ad_{g^{-1}}^{*}\phi(\alpha_p) \backl
(r_{g^{-1}} \pg (g)), ~~ (g^{-1})^* \alpha_p), ~~~ g \in G, ~ p \in P,
\end{equation}
where $x \in \fg, \alpha_p \in T_{p}^{*}P$.
It has the following properties:

(1) The action of $\fg$ on $\Gamma(A)$ induced by this $G$-action on $A$
coincides with the action given by (\ref{eq_fg-on-A}).

(2) The action of $G$ commutes with the anchor map $-\sigma - \ps$ of $A$, and
$G$ acts on $\Gamma(A)$ as automorphisms with respect to the Lie bracket $\{ ~
, ~ \}$.

(3) For any $x \in \fg$ and $g \in G, ~ i(Ad_g x) = g \cdot i(x)$.

In other words, $A$ is a Harish-Chandra $G$-Lie algebroid \cite{bb:jantzen}.

(4) The action of $G$ preserves the scalar product on $A$ given by
\begin{equation}
\label{eq_scalar-on-Ap}
\langle x + \alpha_p, ~ y + \beta_p \rangle_p ~ = ~ (\beta_p, ~\sigma_x(p))
 ~ + ~ (\alpha_p, ~\sigma_y(p)) ~ = ~ \phi(\beta_p)(x) ~ + ~
\phi(\alpha_p)(y).
\end{equation}
\end{prop}

\noindent
{\bf Proof}. The proof is straightforward. We just point out that
to prove that (\ref{eq_G-on-A}) defines an action, one only needs
the fact that the Poisson structure $\pg$ on $G$ makes $G$ into
a Poisson Lie group, i.e.,
\[
\pg(gh) ~ = ~ l_g \pg(h) ~ + ~ r_h \pg(g), ~~~ \forall g, h \in G.
\]
In the proof of (1), one needs the fact that $\ps$ vanishes at
the identity point of $G$. (2) is a reformulation of the fact that
$\sigma$ is a Poisson action. (3) is obvious. (4) is straightforward.
\qed

\begin{prop}
\label{prop_to-fd}

1) the map
\begin{equation}
\label{eq_Phi}
\Phi: ~ A = (P \times \fg) \bowtie T^*P \lrw \fd: ~~ (x, \alpha_p) \Map x +
\phi(\alpha_p)
\end{equation}
is a Lie algebroid morphism from $A$ to the Lie algebra $\fd = \fg \bowtie
\fgs$, considered as a Lie algebroid over a one point space.
With respect to the following action of $G$ on $\fd$ (see (\ref{eq_G-on-fd})
in Section \ref{sec_poi}:
\[
Ad_g (x + \xi) ~ = ~ Ad_g x + Ad_{g^{-1}}^{*} \xi \backl (r_{g^{-1}}
\pg(g)) ~ + ~ Ad_{g^{-1}}^{*} \xi,
\]
where $g \in G, x \in \fg$ and $\xi \in \fgs$, the map $\Phi$ is in fact
a morphism of Harish-Chandra $G$-Lie algebroids;

2) The map $\Phi$ induces the following representation of $A$ on the
trivial vector bundle $P \times \fd$ via the adjoint representation of
$\fd$ on itself:
\begin{equation}
\label{eq_A-on-trivial-fd}
D_{\bar{x}+\alpha}^{A} (\by + \bar{\xi}) ~ = ~ (-\sigma_{\bar{x}} -
\ps \alpha) \cdot (\by + \bar{\xi}) ~ + ~ [\bx + \phi(\alpha), ~ \by
+ \xi]_{\fd}
\end{equation}
where $\by + \bar{\xi}$ is a section of $P \times \fd$
with
$\by \in C^{\infty}(P, \fg)$ and $\bar{\xi} \in C^{\infty}(P, \fgs)$.
This representation of $A$ on $P \times \fd$ preserves the
symmetric product $\la ~,~ \ra_{\fd}$ on $P \times \fd$ given by
(\ref{eq_scalar}) on each fiber.
\end{prop}

\noindent
{\bf Proof.} 1) is easily reduced to the fact that for any $\bx \in C^{\infty}
(P, \fg)$ and any $1$-form $\alpha$ on $P$,
\begin{equation}
\label{eq_pfg-mor}
\phi(D_{\bar{x}} \alpha) ~ = ~ - \sigma_{\bar{x}} \cdot \phi(\alpha) ~ + ~
ad_{\bar{x}}^{*} \phi(\alpha).
\end{equation}
This identity can be proved directly from the definition of $D_{\bar{x}}
\alpha$. If we equip the trivial vector bundle $P \times \fgs$
the representation of the transformation Lie algebroid
$P \times \fg$ given by
\[
D_{\bar{x}} \bar{\xi} ~ = ~ - \sigma_{\bar{x}} \cdot \bar{\xi}
{}~ + ~ ad_{\bar{x}}^{*} ~\bar{\xi}.
\]
Then (\ref{eq_pfg-mor}) is saying that the bundle map
\[
T^*P \lrw P \times \fgs: ~ \alpha_p \Map (p, ~ \phi(\alpha_p))
\]
is a $(P \times \fg)$-morphism. Proof of 2) is straightforward.
\qed

\bigskip
In general, if $A$ is any Lie algebroid over $P$ with anchor map $\rho$,
the kernel $\rho$ in each fiber $A_p$ has an
induced Lie algebra structure, namely, for any $a_p, b_p \in
\ker(\rho|_{A_p})$, extend them arbitrarily to sections $a$ and $b$
of $A$ and define
\[
[a_p, ~ b_p] ~ = ~ \{a, ~ b\}(p).
\]
This is independent of the extensions.

For a transformation Lie algebroid $P \times \fg$, the kernels of the
anchor map $-\sigma$ are the stabilizer Lie subalgebras
\begin{equation}
\label{eq_fgp}
\fg_p ~ = ~ \{ x \in \fg: ~~ \sigma_x (p) = 0\},
\end{equation}
For the cotangent bundle Lie algebroid $T^*P$ of a Poisson manifold $(P, \pi)$,
we get the Lie algebra
\begin{equation}
\label{eq_tp}
{\frak t}_p ~ = ~ \{\alpha_p \in T_{p}^{*}P: ~~ \ps(\alpha_p) = 0\}
\end{equation}
at each point $p \in P$. They are called the transversal Lie algebras
to the symplectic leaves of the Poisson structure on $P$ \cite{we:local}.

For the Lie algebroid $A = (P \times \fg) \bowtie T^*P$ constructed in Theorem
\ref{thm_A}, the kernel of the anchor map $-\sigma-\ps$ in the fiber
$A_p$ over $p \in P$ is the space
\begin{equation}
\label{eq_lp}
{\frak l}_p ~ = ~ \{ (x, \alpha_p): ~ x \in \fg, ~ \alpha_p \in T_{p}^{*}P: ~
\sigma_x(p) + \ps(\alpha_p) = 0 \}.
\end{equation}
It is an
isotropic subspace of $A_p$ with respect to
$\langle ~, ~ \rangle_p$ given by (\ref{eq_scalar-on-Ap}).

\begin{prop}
\label{prop_lp}

(1) For each $p \in P$, the map
\begin{equation}
\label{eq_lp-to-fd}
\Phi: ~ {\frak l}_p \lrw \fd = \fg \bowtie \fgs: ~~ (x, \alpha_p) \Map x +
\phi(\alpha_p)
\end{equation}
is a Lie algebra homomorphism from ${\frak l}_p$ to
the double Lie algebra $\fd$ of $(\fg, \fgs)$.

(2) For any $g \in G$ and $p \in P$, we have $g \cdot {\frak l}_p
= {\frak l}_{gp}$, where $g \cdot {\frak l}_p$ denotes the image of
${\frak l}_p$ in $A_{gp}$ under the action of $g$ on $A$
given by (\ref{eq_G-on-A}).
\end{prop}

The proof of this Proposition is again straightforward.

\begin{prop}
\label{prop_transitive}
When $P$ is symplectic or when $P$ is $G$-homogeneous, the
Lie algebroid $A = (P \times \fg) \bowtie T^*P$
is transitive. Consequently, the subbundle $L$ of $A$
defined by the kernel of the anchor map $-\sigma - \ps$
is a Lie algebra bundle of rank $= \dim \fg$.
\end{prop}

In general, the dimensions of the Lie algebra ${\frak l}_p$'s may
vary. Consider again the
stabilizer Lie subalgebra $\fg_p$ and the transversal Lie
algebra ${\frak t}_p$ at each point $p \in P$.
Clearly,
\[
\fg_p \subset {\frak l}_p, ~~~~~ {\frak t}_p \subset {\frak l}_p
\]
as Lie subalgebras,
and
\[
\fg_p \cap {\frak t}_p ~ = ~ 0.
\]
Set
\begin{eqnarray}
\label{eq_Op}
O_p & = & \{\sigma_x(p): ~ x \in \fg \} \subset T_pP, \\
\label{eq_Sp}
S_p & = & \{\ps (\alpha_p): ~ \alpha_p \in T^{*}_{p}P\} \subset T_pP.
\end{eqnarray}
Then
\beqa
\dim {\frak l}_p & = & \dim \fg + \dim P - \dim (O_p + S_p)\\
& = & (\dim \fg - \dim O_p) + (\dim P - \dim S_p) + \dim(O_p \cap S_p)\\
& = & \dim \fg_p ~ + ~ \dim {\frak t}_p ~ + ~
\dim(O_p \cap S_p).
\eeqa

\begin{exam}
\label{exam_symplectic-iso}
{\em
When $P$ is symplectic, i.e., when $\ps: T_{p}^{*}P \rightarrow T_pP$ is
one-to-one for each $p \in P$, the map
$\Phi: {\frak l}_p \rightarrow \fd$ is injective, so $\Phi({\frak l}_p)$
is a maximal isotopic Lie subalgebra of $\fd$.

Let $\omega$ be the symplectic $2$-form on $P$ related to $\pi$ by
\[
 \omega(\ps \alpha, ~ \ps \beta) ~ = ~ \pi (\alpha, ~ \beta).
\]
Then
\[
\Phi({\frak l}_p) ~ = ~ \{x - x\backl \Phi(\omega_p): ~~ x \in \fg \} \subset
\fd,
\]
where $\Phi(\omega_p) = (\Phi \wedge \Phi)(\omega_p)
\subset \fgs \wedge \fgs$. Notice that $\Phi({\frak l}_p)$
is transversal to $\fgs$ in $\fd$. Thus $(\fd, \fgs, \Phi({\frak l}_p))$
form a Manin triple \cite{dr:bigbra}, and $(\fgs, \Phi({\frak l}_p))$ becomes
a Lie bialgebra.

Let $(G^*, \pgs)$ be the simply-connected Poisson Lie group
dual to $G$. The tangent Lie bialgebra of $G^*$ is
$(\fgs, \fg)$.
Define
\[
\pi_{\mbox{\tiny G}^{*}, p} ~ = ~ \pi_{\mbox{\tiny G}^{*}} ~ + ~
(\Phi(\omega_p))^r ~ - ~ (\Phi(\omega_p))^l,
\]
where $(\Phi(\omega_p))^r$ (resp. $(\Phi(\omega_p))^l$) is the right invariant
bi-vector field on $G^*$ with value $\Phi(\omega_p) \in \fgs \wedge \fgs$
at the identity element $e \in G^*$. Then $(G^*, \pi_{\mbox{\tiny G}^{*}, p})$
is also a Poisson Lie group, and its tangent Lie bialgebra
is $(\fgs, \Phi({\frak l}_p))$ \cite{lu:thesis} \cite{daso:affine}.

The bi-vector field $\pi_{\mbox{\tiny G}^{*}} + (\Phi(\omega_p))^r$ is also
Poisson. It is an example of an affine Poisson
structure on $G^*$ \cite{lu:thesis} \cite{daso:affine}. When $P$ is
connected and simply-connected, it is shown in \cite{lu:thesis}
that for any $p \in P$ there exists a unique map $J_p: P \rightarrow G^*$
such that

(1) $J_p(p) = e$;

(2) $\ps (J_{p}^{*} x^l) = \sigma_x$ for every $x \in \fg$,
where $x^l$ is the left invariant $1$-form on $G^*$ such that
$x^l (e) = x$;

(3) $J_p: (P, \pi) \rightarrow (G^*, \pi_{\mbox{\tiny G}^{*}} +
(\Phi(\omega_p))^r)$ is a Poisson map.

The map $J_p$ is called a moment map for the Poisson action $\sigma$ of
$G$ on $P$ \cite{lu:mom}.

If $p_1$ is any other point in $P$, the maps $J_p$ and $J_{p_1}$
are related by
\[
J_p (q) ~ = ~ J_p (p_1) ~ J_{p_1} (q), ~~~ q \in P,
\]
and
\[
\Phi(\omega_p) ~ = ~ -r_{u^{-1}} \pi_{\mbox{\tiny G}^{*}} (u) ~ + ~ Ad_u \Phi(
\omega_{p_1}).
\]
The Manin triples $(\fd, ~ \fgs, ~ \Phi({\frak l}_p))$ and
$(\fd, ~ \fgs, ~ \Phi({\frak l}_{p_1}))$ are related by
\[
(\fd, ~ \fgs, ~ \Phi({\frak l}_p)) ~ = ~
Ad_{J_p(p_1)} \left(\fd, ~ \fgs, ~ \Phi({\frak l}_{p_1})\right)
\]
where $J_p(p_1) \in G^*$ acts on $\fd$ by (\ref{eq_Gs-on-fd}).
}
\end{exam}

\begin{exam}
\label{exam_dressing-orbit}
{\em
The left dressing vector fields on $G^*$
are defined by
\begin{equation}
\label{eq_left-dressing}
\sigma_x (u) ~ = ~ \pi_{G^*}^{\#} (x^l)
\end{equation}
where $x \in \fg$.
The map
\[
\fg \lrw \chi(G^*): ~ x \Map \sigma_x
\]
defines a Lie algebra anti-homomorphism from
$\fg$ to the Lie algebra $\chi(G^*)$ of vector fields
on $G^*$. Assume that these vector fields integrate to
an action of $G$ on $G^*$. It is called the left dressing
action of $G$ on $G^*$ \cite{sts:dressing}. It is a
Poisson action. It in fact has the identity map $G^* \rightarrow
G^*$ as a moment map \cite{lu:mom}. The orbits of
the dressing action are exactly the symplectic leaves of
$G^*$. Let $P$ be such an orbit. Then the
action $G \times P \rightarrow P$ of $G$ on $P$ is a Poisson action.

For $u \in P \subset G^*$, the Lie algebra ${\frak l}_u$ can
be identified with the Lie subalgebra $Ad_{u^{-1}} \fg$ of $\fd$, where, again,
$Ad_{u^{-1}}: \fd \rightarrow \fd$ denotes the action of $ u^{-1} \in
G^*$ on $\fd$ given by (\ref{eq_Gs-on-fd}).
}
\end{exam}

\begin{exam}
\label{exam_homog-1}
{\em
When $P$ is a homogeneous $G$-space, we have $\dim {\frak l}_p
= \dim \fg$ for each $p \in P$, and the map $\Phi:
{\frak l}_p \rightarrow \fd$ is injective. Therefore,
each $\Phi({\frak l}_p)$ is a maximal isotropic Lie subalgebra of $\fd$.
They are used in \cite{dr:homog} to
classify Poisson homogeneous $G$-spaces. We will treat them in
more details in Section \ref{sec_homog}.
}
\end{exam}

\begin{exam}
\label{exam_pi-0}
{\em
When the Poisson bi-vector field $\pi$ on $P$ vanishes at
a certain point $p \in P$, we have
\[
{\frak l}_p ~ = ~ \{(x, ~ \alpha_p): ~ x \in \fg, ~ \alpha_p \in T_{p}^{*}P,
{}~ \sigma_x (p) = 0 \} ~\cong ~ \fg_p \oplus T_{p}^{*}P.
\]
To describe the Lie algebra structure on ${\frak l}_p$, we first note
that the action of $G_p$ on $P$, where $G_p$ is the stabilizer subgroup
of $G$ at $p$, linearizes to an action of $G_p$ on $T_{p}P$
and thus on $T_{p}^{*}P$. Denote the
corresponding action of $\fg_p$ on $T_{p}^{*}P$ by
\begin{equation}
\label{eq_fgp-on-tsp}
\fg_p \times T_{p}^{*} P \lrw T_{p}^{*}P: ~ (x, ~ \alpha_p) \Map
x \cdot \alpha_p.
\end{equation}
On the other hand, since the action is Poisson, we know
from
(\ref{eq_poisson-action-infinitesimal})
that
\[
\sigma_{\delta(x)} (p) ~ = ~ 0, ~~~~~\forall x \in \fg_p.
\]
Thus
\[
(x, ~ [\phi(\alpha_p), ~ \phi(\beta_p) ]_{\fgs} ) ~ = ~ 0,
{}~~~~~~~\forall x \in \fg_p, ~ \alpha_p, ~ \beta_p \in T_{p}^{*}P.
\]
In other words,
\[
ad_{\phi(\alpha_p)}^{*} \fg_p \subset \fg_p, ~~~~~~~ \forall \alpha_p
\in T_{p}^{*}P.
\]
The Lie algebra structure on ${\frak l}_p$ is now given by
\begin{equation}
\label{eq_bra-on-lp}
[(x, \alpha_p), ~ (y, \beta_p)] ~ = ~ ([x, y] + ad_{\phi(\alpha_p)}^{*} y
- ad_{\phi(\beta_p)}^{*} x, ~~ [\alpha_p, \beta_p] + x \cdot \beta_p
- y \cdot \alpha_p).
\end{equation}

Note that the Lie algebra ${\frak l}_p = \fg_p + T_{p}^*P$
is an example of a double Lie algebra \cite{lu-we:poi} \cite{ks:podr}
(called a matched pair of Lie algebras in \cite{mj:mat})
in the sense that it contains $\fg_p$ and $T_{p}^*P$ as Lie subalgebras
and it is isomorphic to $\fg_p \oplus T_{p}^*P$ as vector spaces.
We can think of ${\frak l}_p$ as being built out of the
Lie algebras $\fg_p$ and $T_{p}^*P$ together with the action of $\fg_p$ on
$T_{p}^{*}P$ given by (\ref{eq_fgp-on-tsp}) and the action of
$T_{p}^{*}P$ on $\fg_p$ by
\[
T_{p}^{*} P \times \fg_p \lrw \fg_p: ~~ (\alpha_p, ~ x)
\Map ad_{\phi(\alpha_p)}^{*} x.
\]
}
\end{exam}

\begin{exam}
\label{exam_coiso-1}
{\em
Combining the situations in Example \ref{exam_homog-1} and
Example \ref{exam_pi-0}, we consider the case
when $P$ is a Poisson homogeneous $G$ space and the
Poisson bivector field $\pi$ on $P$
vanishes at some point $p \in P$. In this case, we
can identify $P$ with the quotient space $G/G_p$,
where $G_p$ is the stabilizer subgroup of $G$ at $p$. The Poisson
structure $\pi$ is the unique one on $G/G_p$ such that the projection
from $G$ to $G/G_p$ is a Poisson map. From Example \ref{exam_pi-0}
we see that
\[
\fg_{p}^{\perp} ~ = ~ \{\xi \in \fgs: ~~ (\xi, x ) = 0 ~~\forall x \in \fg_p \}
\subset \fgs
\]
is a Lie subalgebra of $\fgs$. The Lie subalgebra $\fg_p$ of $\fg$ is said
to be coisotropic relative to the Lie bialgebra $(\fg, \fgs)$ \cite{lu-we:poi}.
The Lie algebra ${\frak l}_p$ is now isomorphic to the Lie subalgebra
$\fg_p + \fg_{p}^{\perp}$ of $\fd = \fg \bowtie \fgs$.

The fact that
the Lie algebra $ {\frak l}_p =
\fg_p + \fg_{p}^{\perp}$  is
a double Lie algebra
will be used in Section \ref{sec_homog} to
calculate the relative Lie algebra cohomology of ${\frak l}_p$
relative to the Lie subalgebra $\fg_p$, which will be shown to
be isomorphic to the $G$-invariant Poisson cohomology of the
Poisson homogeneous space $G/G_p$. See
Theorem \ref{thm_coisotropic-cohomo}.

A special case of this situation is when $\fg_{p}^{\perp}$ is a Lie ideal of
$\fgs$. In such a case, $G_p$ is a Poisson Lie subgroup of $G$. The Lie
bracket on $\fg_p + \fg_{p}^{\perp}$ is a little simpler:
\begin{equation}
\label{eq_poi-sub}
[x + \xi, ~ y + \eta] ~ = ~ [x, ~y] ~ + ~ [\xi, ~ \eta] ~ + ~ ad_{x}^{*} \eta
{}~ - ~ ad_{y}^{*} \xi
\end{equation}
for $x, y \in \fg_p$ and $\xi, \eta \in \fg_{p}^{\perp}$. In particular,
$\fg_{p}^{\perp}$ is a Lie ideal of $ {\frak l}_p \cong \fg_p +
\fg_{p}^{\perp}$.

This example will be carried out further in
Example \ref{exam_coiso-2} in Section \ref{sec_homog}.
}
\end{exam}

\section{Examples}
\label{sec_exam}

In this section, we give examples of the Lie algebroid $A = (P \times \fg)
\bowtie T^*P$ constructed in Section \ref{sec_A}.

\begin{exam}
\label{exam_0}
{\em
Any Lie group $G$ can be regarded as a Poisson Lie group with
the zero Poisson structure. In this case, an action of
$G$ on a Poisson manifold is Poisson if and only if each $g \in G$
acts on $P$ as a Poisson diffeomorphism. In particular,
each $x \in \fg$ acts on the space of $1$-forms on $P$ as
a derivation. The corresponding Lie algebroid structure
on the vector bundle $(P \times \fg) \oplus T^*P$ described in
Theorem \ref{thm_A} becomes a
semi-direct product Lie algebroid structure \cite{h-mk:cat}.
}
\end{exam}

\begin{exam}
\label{exam_G}
{\em
Let $G$ be a Poisson Lie group.
Consider the left action of $G$ on itself by left translations.
This is a Poisson action by definition. By trivializing
the cotangent bundle $T^*G$ to $G \times \fgs$ via,
\[
T^*G \lrw G \times \fgs: ~~ \xi_g \Map  \phi(\xi_g) ~ = ~ r_{g}^{*} \xi_g, ~~~~
g \in G, ~
\xi_g \in T_{g}^{*}G,
\]
we can trivialize the vector bundle $(P \times \fg) \oplus T^*G$ to
$G \times (\fg \oplus \fgs) = G \times \fd$ by the
map
\[
\Phi: ~ (x, ~ \xi_g) \Map x ~ + ~ \phi(\xi_g) ~ = ~ x ~ + ~ r_{g}^{*} \xi_g,
{}~~~~ g \in G, ~ \xi_g \in T_{g}^{*}G.
\]
In this case, the
Lie algebroid structure on the bundle $A$ described in Theorem
\ref{thm_A} is the transformation Lie algebroid structure on $G \times
\fd$ given by the left infinitesimal action of $\fd$ on $G$
given by (\ref{eq_fd-on-G-left}).

The maximal isotropic Lie subalgebra ${\frak l}_g$
of $\fd$ in this case is simply the Lie subalgebra
$Ad_g \fgs$, the image of $\fgs$ under the action of $g$
on $\fd$ given by (\ref{eq_G-on-fd}).
}
\end{exam}

\begin{exam}
\label{exam_homog-2}
{\em
Assume that the action $\sigma: G \times P \rightarrow P$
is transitive, so $P$ is a Poisson homogeneous $G$-space. In this case,
the map
\[
T^*P \lrw P \times \fgs: ~~ \alpha_p \Map (p, ~ \phi(\alpha_p))
\]
embeds $T^*P$ into $P \times \fgs$ as a subbundle whose fiber over
the point $p \in P$ is
\[
\fg_{p}^{\perp} ~ = ~ \{ \xi \in \fgs: ~~ (\xi, x) = 0, ~ \forall
x \in \fg_p \}.
\]
Thus via the map
\[
\Phi: ~ A = (P \times \fg) \bowtie T^*P \lrw \fd ~ = ~
\fg \bowtie \fgs: ~~ (x, \alpha_p) \Map x + \phi
(\alpha_p)
\]
we can also embed $A$ into
$P \times (\fg \oplus \fgs)$
 as a subbundle
whose fiber over $p \in P$ is $\fg \oplus \fg_{p}^{\perp}$.
Under such an embedding, a section of $A$ has the form $a = \bx + \bar{\xi}$
where $\bx \in C^{\infty}(P, \fg)$ and $\bar{\xi} \in C^{\infty}(P, \fgs)$
is such that $\bar{\xi}(p) \in \fg_{p}^{\perp}$ for all $p \in P$. Denote by
$\ps(\bar{\xi})$ the vector field $\ps(\alpha)$ if
$\phi(\alpha) = \bar{\xi}$. The
anchor of $A$ is then the vector field
$-\sigma_{\bar{x}} - \ps(\bar{\xi})$.
Since $\Phi$ is a Lie algebroid morphism from $A$ to the
Lie algebra $\fd = \fg \bowtie \fgs$ by Proposition \ref{prop_to-fd},
the Lie bracket on the sections of
$A$ now looks like that for a
transformation Lie algebroid:
\begin{equation}
\label{eq_bra-on-A-homog}
\{\bx + \bar{\xi}, ~ \by + \bar{\eta}\} ~ = ~
[\bx + \bar{\xi}, ~ \by + \bar{\eta}]_{\frak d} ~ - ~
(\sigma_{\bar{x}} + \ps(\bar{\xi})) \cdot (\by + \bar{\eta}) ~ + ~
(\sigma_{\bar{y}} + \ps (\bar{\eta})) \cdot (\bx + \bar{\xi}).
\end{equation}
}
\end{exam}

\begin{exam}
\label{exam_dressing}
{\em
Consider now the example of the left dressing action of $G$ on its dual
Poisson Lie group $(G^*, \pgs)$ (see Example \ref{exam_dressing-orbit}).

Identify $T^{*}G^*$ with the trivial vector bundle
$G^* \times \fg$ by left translations on $G^*$. Then the
cotangent bundle Lie algebroid $T^{*}G^*$ becomes the same as the
transformation Lie algebroid $G^* \times \fg$ defined by the
left dressing action of $\fg$ on $G^*$. Thus the Lie algebroid
$A$, whose underlying vector bundle is now identified
with the trivial vector bundle $G^* \times (\fg \oplus \fg)$,
is built out of two copies of the transformation Lie algebroid
$G^* \times \fg$ and a pair of representations $D$ and $D'$ of it on itself.
These two representations are respectively given by
\begin{eqnarray}
\label{eq_pair-1}
D_{\bar{x}} \bar{y} & = & -\sigma_{\bar{x}} \cdot \bar{y} ~ + ~
[\bx, ~ \by]_{\fg} ~ + ~ ad_{\tau(\bar{y})}^{*} \bx \\
\label{eq_pair-2}
D_{\bar{y}}^{'} \bx & = & - \sigma_{\bar{y}} \cdot \bx ~ + ~
 ad_{\tau(\bar{y})}^{*} \bx,
\end{eqnarray}
where
\begin{equation}
\label{eq_tau}
\tau: ~ C^{\infty}(G^*, \fg) \lrw C^{\infty}(G^* \fgs): ~~
\tau(\by) (u) ~ = ~ \by(u) \backl (l_{u^{-1}}  \pgs (u)), ~~~ u \in G^*.
\end{equation}
The resulting Lie algeroid structure on $G^* \times (\fg \oplus \fg)$
has the property that the bundle map
\[
G^* \times (\fg \oplus \fg) \lrw G^* \times \fg: ~~ (x_u, ~ y_u)
\Map x_u ~ + ~ y_u
\]
is a Lie algebroid homomorphism.
}
\end{exam}

\section{Proof of Theorem}
\label{sec_proof}

We give the proof of Theorem \ref{thm_A} in this section.
Following the definition of a Lie algebroid, we need to prove

1) for any $f \in C^{\infty}(P), ~ \bx \in C^{\infty}(P, \fg)$
and any $1$-form $\alpha$ on $P$,
\beqa
\{f \bx, ~ \alpha\} & = & f \{\bx, ~ \alpha\} ~ + ~ (\ps(\alpha) f) \bx\\
\{\bx, ~ f \alpha\} & = & f \{\bx, ~ \alpha\} - (\sigma_{\bar{x}} f) \alpha;
\eeqa

2) the bracket on the sections of $A$ described in Theorem \ref{thm_A}
satisfies the Jacobi identity;

3) the map $-\sigma-\ps$ is a Lie algebra homomorphism from the Lie algebra
of sections of $A$ to the Lie algebra of vector fields on $P$.

\bigskip
1) follows from the fact that $D_{\alpha} \bx$
and $D_{\bar{x}} \alpha$ define representations of Lie algebroids.
This fact can also be used to reduce 2) to Lemma \ref{lem_D-bra}
and 3) to Lemma \ref{lem_anchor-A} as follows.

\begin{lem}
\label{lem_D-bra}
For any $\bx, \by \in C^{\infty}(P, \fg)$ and any $1$-forms $\alpha$ and
$\beta$
on $P$, we have
\begin{eqnarray}
\label{eq_D1}
& & D_{\alpha} \{\bx, ~ \by\} ~ = ~ \{D_{\alpha} \bx, ~ \by\} ~ + ~ \{\bx, ~
D_{\alpha} \by\} ~ + ~ D_{D_{\by} \alpha} \bx ~ - ~
D_{D_{\bx} \alpha} \by \\
\label{eq_D2}
& & D_{\bar{x}} \{\alpha, ~ \beta\} ~ = ~ \{D_{\bar{x}} \alpha, ~ \beta\}
{}~ + ~ \{\alpha, ~ D_{\bar{x}} \beta\} ~ + ~ D_{D_{\beta} \bar{x}} \alpha
{}~ - ~ D_{D_{\alpha} \bar{x}} \beta.
\end{eqnarray}
\end{lem}

\begin{lem}
\label{lem_anchor-A}
For any $\bx \in C^{\infty}(P, \fg)$ and any $1$-form $\alpha$ on $P$,
we have
\begin{equation}
\label{eq_anchor-A}
\sigma_{D_{\alpha} \bar{x}} ~ - ~ \ps(D_{\bar{x}} \alpha) ~
= ~ [\sigma_{\bar{x}}, ~ \ps(\alpha)].
\end{equation}
\end{lem}

We prove Lemma \ref{lem_anchor-A} first, since it is easier to
prove and it will be used in the proof of Lemma \ref{lem_D-bra}.

\bigskip
\noindent
{\bf Proof of Lemma \ref{lem_anchor-A}} Fix the $1$-form $\alpha$ on
$P$. For any $\bx \in C^{\infty}(P, \fg)$, set
\[
V_{\alpha} (\bx) ~ = ~ \sigma_{D_{\alpha} \bar{x}} ~ - ~ \ps(D_{\bar{x}}
\alpha) ~ - ~ [\sigma_{\bar{x}}, ~ \ps(\alpha)].
\]
It follows from the axioms for Lie algebroid representations that
\[
V_{\alpha} (f \bar{x}) ~ = ~ f V_{\alpha} (\bar{x})
\]
for any $f \in C^{\infty}(P)$. Thus it suffices to show that
$V_{\alpha}(x) = 0$ for any constant function $x \in C^{\infty}(P, \fg)$
corresponding to $x \in \fg$. In other words, we need to show that
\begin{equation}
\label{eq_const}
\sigma_{ad_{\phi(\alpha)}^{*} x} ~ + ~ \ps (L_{\sigma_{x}} \alpha) ~ - ~
[\sigma_{x}, ~ \ps (\alpha) ] ~ = ~ 0
\end{equation}
for any $x \in \fg$. By pairing the left hand side with an arbitrary
$1$-form $\beta$ on $P$, we see that this is equivalent to
\[
(L_{\sigma_x} \pi) (\alpha, ~ \beta) ~ = ~ (x, ~ [\phi(\alpha), ~
\phi(\beta)]_{{\frak g}^*}),
\]
which is just the infinitesimal condition
(\ref{eq_poisson-action-infinitesimal})
for the action $\sigma$ to be Poisson.
\qed

\noindent
{\bf Proof of Lemma \ref{lem_D-bra}}. For $\bx, ~\by, ~\alpha$ and $\beta$ as
given, set
\beqa
& & V_{\alpha}(\bx,\by) ~ = ~ D_{\alpha} \{\bx, ~ \by\}
{}~ - ~ \{D_{\alpha} \bx, ~ \by\} ~ - ~ \{\bx, ~
D_{\alpha} \by\} ~ - ~ D_{D_{\by}\alpha} \bx ~ + ~
D_{D_{\bx}\alpha} \by \\
& & W_{\alpha, \beta} (\bx) ~ = ~ D_{\bar{x}} \{\alpha, ~ \beta\}
{}~ - ~ \{D_{\bar{x}} \alpha, ~ \beta\}
{}~ - ~ \{\alpha, ~ D_{\bar{x}} \beta\} ~ - ~ D_{D_{\beta} \bar{x}} \alpha
{}~ + ~ D_{D_{\alpha} \bar{x}} \beta.
\eeqa
Using Identity (\ref{eq_const}) and the axioms for Lie algebroid
representations, we see that
\[
V_{\alpha}(f \bx, ~ \by) ~ = ~V_{\alpha}(\bx, ~ f\by) ~ = ~
f V_{\alpha}(\bx, ~ \by)
\]
for any $f \in C^{\infty}(P)$. Thus it suffices to show that
$V_{\alpha}(x, y) = 0$ for any constant functions $x, y \in
C^{\infty}(P, \fg)$ corresponding to $x, y \in \fg$. But this can be deduced
from the following fact about the Lie bialgebra $(\fg, \fgs)$:
\begin{equation}
\label{eq_for}
ad_{\xi}^* [x, y] ~ = ~ [ad_{\xi}^{*} x, ~ y] ~ + ~ [x, ~ ad_{\xi}^{*} y]
{}~ + ~ ad_{ad_{y}^{*} \xi}^{*}x ~ - ~ ad_{ad_{x}^{*} \xi}^{*} y
\end{equation}
for any $x, y \in \fg$ and $\xi \in \fgs$. The proof of (\ref{eq_for})
can be obtained from the explicit formula (\ref{eq_bra-on-d})
for the Lie bracket on the double $\fd$ of $(\fg, \fgs)$. It
can also be found in \cite{am:thesis}.

For $W_{\alpha, \beta}$, it again follows from the axioms for
Lie algebroid representations that
\[
W_{\alpha, \beta} (f \bx) ~ = ~ f W_{\alpha, \beta} (\bx)
\]
for any $f \in C^{\infty}(P)$. Thus it suffices to show that
$W_{\alpha, \beta}(x) = 0$ for constant functions $x \in C^{\infty}(P, \fg)$
corresponding to $x \in \fg$. In other words, we need to show that
\[
L_{\sigma_x} \{ \alpha, ~ \beta\} ~ - ~ \{L_{\sigma_x} \alpha, ~ \beta\}
{}~ - ~ \{\alpha, ~ L_{\sigma_x} \beta\} ~ + ~ D_{\tilde{x}_{\beta}} \alpha
{}~ - ~ D_{\tilde{x}_{\alpha}} \beta ~ = ~ 0,
\]
where
\[
\tilde{x}_{\beta} ~= ~ ad_{\phi(\beta)}^{*} x, ~~~~~~~~
\tilde{x}_{\alpha} ~= ~ ad_{\phi(\alpha)}^{*} x.
\]
Using Formula (\ref{eq_2}) for the
Lie bracket on $1$-forms on $P$ and Identity (\ref{eq_const}), we
can reduce $W_{\alpha, \beta}(x) = 0$ to
\[
(L_{\sigma_x} \pi) (\alpha, ~ \beta) ~ = ~ (x, ~ [\phi(\alpha),
{}~ \phi(\beta) ]_{\fgs} )
\]
which is nothing but the infinitesimal criterion
(\ref{eq_poisson-action-infinitesimal}) for $\sigma$ being
Poisson.
\qed

\section{Poisson homogeneous spaces}
\label{sec_homog}

In this section, we assume that $G$ is a connected Poisson Lie group
and that $P$ is a homogeneous $G$-space with a Poisson structure
$\pi$ such that the action of $G$ on $P$ is  Poisson.
We will treat two aspects of such a Poisson manifold:
its symplectic leaves and its
$G$-invariant Poisson cohomology.  We first look at its symplectic leaves.

For each $p \in P$, embed ${\frak l}_p$ into $\fd$ via
\[
\Phi: ~ {\frak l}_p \lrw \fd: ~ (x, \alpha_p) \Map x + \phi(\alpha).
\]
Then we get a Lie algebra homomorphism from ${\frak l}_p$ to $\chi(G)$:
\[
{\frak l}_p \lrw \chi(G): ~ (x, \alpha_p) \Map
\lambda_{x+\phi(\alpha_p)}(g) ~ = ~ r_g \left(
Ad_g x + (Ad_{g^{-1}}^{*} \phi(\alpha_p) \backl (r_{g^{-1}}
\pg(g)) \right).
\]
Recall (see (\ref{eq_fd-on-G-right}) in Section \ref{sec_poi}) that
the vector fields $\lambda_{x + \xi}$, for $x \in \fg, ~
\xi \in \fgs$,
define a right action of $\fd$ on $G$.
Set
\begin{equation}
\label{eq_sigma-p}
\sigma_p: ~ G \lrw P: ~~ g \Map gp.
\end{equation}

\begin{lem}
\label{lem_distribution} For any $p \in P, ~ g \in G$,
the image of the subspace
\[
\{\lambda_{x+\phi(\alpha_p)}(g): ~~ (x, \alpha_p) \in {\frak l}_p \}
\subset T_g G
\]
of $T_g G$ under the differential of the
map $\sigma_p$ is exactly the subspace
$im (\ps|_{T_{gp}^{*}P})$ of $T_{gp} P$.
\end{lem}

\noindent
{\bf Proof}. We know from Proposition \ref{prop_harish-chandra-A} and
Proposition \ref{prop_lp} that
\[
\sigma_p (\lambda_{x+\phi(\alpha_p)}(g)) ~ =~ - \ps((g^{-1})^* \alpha_p)
\in im (\ps|_{T_{gp}^{*}P}).
\]
The fact that we get all elements in $im (\ps|_{T_{gp}^{*}P})$ this way is
because the action of $G$ on $P$ is transitive.
\qed

We now have the following description of the symplectic leaves in $P$.

\begin{thm}
\label{thm_leaves}
Fix any $p \in P$. Assume that $L_p$ is a connected Lie group
and that the vector fields
\[
\lambda_{x+\phi(\alpha_p)}, ~~~ (x, \alpha_p) \in {\frak l}_p
\]
integrate to a right action of $L_p$ on $G$. Then the images of
the $L_p$ orbits in $G$ under the map
\[
\sigma_p: ~ G \lrw P: ~~ g \Map gp
\]
give all the symplectic leaves in $P$.
\end{thm}

{\bf Proof.} By definition, the symplectic leaves of $P$ are the integral
submanifolds of the distribution $\{im (\ps|_{T_{p}^{*}P})\}_{p \in P}$
on $P$. By Lemma \ref{lem_distribution}, we know that this distribution is the
image under $\sigma_p$ of the distribution on $G$ defined by
the vector fields $\lambda_{x+\phi(\alpha_p)}, (x, \alpha_p) \in {\frak l}_p,$
on $G$. The integral submanifolds of the latter distribution
are exactly all the $L_p$-orbits on $G$. Thus these orbits
project to $P$ under $\sigma_p$ to give all the symplectic leaves in $P$.
\qed

\begin{rem}
\label{rem_after-leaves}
{\em
Note that when the double group $D$ of $G$ and $G^*$ have
a global decomposition $D = G G^*$, the right action of $L_p$ on $G$
described in Theorem \ref{thm_leaves} is simply the restriction to
$L_p$ of
the right action of $D$ on $G = G^* \backslash D$ (see (\ref{eq_D-on-G-right})
in Section \ref{sec_poi}).
}
\end{rem}

\bigskip
We now turn to the $G$-invariant Poisson cohomology of $P$.
Recall that the Poisson cohomology of $P$ is by definition the
Lie algebroid cohomology of the cotangent bundle Lie algebroid $T^*P$
with trivial coefficients. The corresponding coboundary operator
on $\wedge^{\bullet} TP$ is also given
by
\[
V \Map [\pi, ~ V]
\]
where $[ ~~, ~~]$ denotes the Schouten bracket \cite{ku:schouten}. If $V$ is
$G$-invariant, then so is $[\pi, V]$. This is because
\[
L_{\sigma_x} [\pi, ~ V] ~ = ~ [L_{\sigma_x}, ~ V] ~ = ~ [\sigma_{\delta(x)},
{}~ V]
{}~ = ~ 0
\]
for any $x \in \fg$ (see (\ref{eq_poisson-action-infinitesimal})).
Hence, since $G$ is connected, the space $\Gamma(\wedge TP)^{G}$ of
$G$-invariant multi-vector fields on $P$
is closed under the coboundary operator $[\pi, ~ \bullet]$.

\begin{dfn}
\label{dfn_invariant-cohom}
{\em
The cohomology of $(\Gamma(\wedge TP)^{G}, ~ [\pi, ~ \bullet])$
is called the $G$-invariant Poisson cohomology of $P$. We denote
it by $H_{\pi, G}(P)$.
}
\end{dfn}

\bigskip
Assume now that $P$ is a Poisson homogeneous $G$-space.
We have seen in Example \ref{exam_homog-1} that the
kernel of the anchor map
$-\sigma-\ps$ at each $p \in P$ is a Lie algebra
${\frak l}_p$ of dimension $n = \dim \fg$.
Let $G_p$ be the stabilizer Lie subgroup of $G$
at $p$, and let $\fg_p$ be its Lie algebra. The group
$G_p$ acts on  ${\frak l}_p$ via the action of $G$ on
$A$ given by (\ref{eq_G-on-A}), and the corresponding action
of $\fg_p$ is the adjoint action of $\fg_p$ on ${\frak l}_p$
via the Lie algebra embedding $\fg_p \subset {\frak l}_p$.
Denote by $H^{\bullet}({\frak l}_p, G_p)$ the relative Lie algebra
cohomology of ${\frak l}_p$ relative to $G_p$.

\bigskip
\begin{thm}
\label{thm_invariant-relative}
\[
H_{\pi, G}^{\bullet}(P) \cong H^{\bullet}({\frak l}_p, G_p)
\]
for every $p \in P$.
\end{thm}

\bigskip
\noindent
{\bf Proof.} Denote by $L$
the kernel of $-\sigma-\pi$, so $L$ is a subbundle of $A$
whose fiber at $p \in P$ is the Lie algebra
${\frak l}_p$. By taking fiber-wise Lie bracket of sections of $L$,
we can think of $L$ as a Lie algebroid itself with
the zero anchor map. As a such, it is a ``totally intransitive" Lie
algebroid \cite{mk:book}. The inclusion of $L$ into $A$ makes it
into a Lie subalgebroid of $A$. We know from
Theorem \ref{thm_A} and Proposition \ref{prop_lp} that $L$ is a
$G$-vector bundle and that the action of $G$ on $L$ preserves
the Lie brackets on the fibers of $L$.

Consider now the coboundary operator
$d_{\scriptscriptstyle L}$ for the Lie algebroid cohomology
of $L$ with trivial coefficients (see (\ref{eq_dk})). Since
$L$ has the zero anchor map, the operator $d_{\scriptscriptstyle L}$
is $C^{\infty}(P)$-linear, i.e., $d_{\scriptscriptstyle L}$
is given by the field $d_{{\frak l}_p}, p \in P$, of
fiber-wise operators, where $d_{{\frak l}_p}$ is the
coboundary operator for the cohomology of the Lie algebra
${\frak l}_p$.

For each $p \in P$, set
\begin{equation}
\label{eq_g0}
{\frak g}_{p}^{0} ~ = ~ \{f \in {\frak l}_{p}^{*}: ~ f(x) = 0, ~ \forall
x \in \fg_p \} \subset {\frak l}_{p}^{*}.
\end{equation}
Then the subspace $(\wedge^{\bullet} {\frak g}_{p}^{0})^{G_p}$ of
$G_p$-invariant vectors in $\wedge^{\bullet} {\frak g}_{p}^{0}$
is invariant under $d_{{\frak l}_p}$. The relative Lie algebra
cohomology of ${\frak l}_p$ relative to $G_p$ is by definition
the cohomology of the cochain complex $((\wedge^{\bullet} {\frak
g}_{p}^{0})^{G_p}, d_{{\frak l}_p})$.

Denote by $L_{\sigma}^{0}$ the subbundle of $L^*$ whose fiber at $p$
is ${\frak g}_{p}^{0}$. It is a $G$-invariant subbundle of $L^*$.
By dimension counting,
the map
\[
TP \lrw L^*: ~~ v_p \Map f_{v_p}: ~ (x, ~ \alpha_p) \Map (v_p, ~ \alpha_p)
\]
gives a bundle isomorphism from $TP$ to $L_{\sigma}^{0}$. It is
also a $G$-bundle morphism. Therefore we get a vector space isomorphism
\begin{equation}
\label{eq_theta}
\Gamma(\wedge^{k} TP)^G \lrw \Gamma(\wedge^k L_{\sigma}^{0})^G
\subset \Gamma(\wedge^k L^*),
\end{equation}
where both sides denote $G$-invariant sections in the corresponding
vector bundles.
Denote this map by $\theta$. It remains to show that
\[
d_{\scriptscriptstyle L} \theta(V) ~ = ~ \theta([\pi, V])
\]
for any $V \in \Gamma(\wedge^{k} TP)^G $, for then $\theta$
would give a cochain complex isomorphism from $\Gamma(\wedge^{k} TP)^G $
to $((\wedge^{\bullet} {\frak g}_{p}^{0})^{G_
p}, d_{{\frak l}_p})$ for each $p \in P$, and thus inducing
an isomorphism between their cohomology spaces.

Let $V$ be a $G$-invariant $k$-vector field on $P$. Then
\[
D_{\bar{x}} V ~ = ~ 0, ~~~ \forall \bx \in C^{\infty}(P, \fg).
\]
Consequently, for any $1$-forms $\alpha_i, i = 1, ..., k, $ on $P$,
we have
\[
-\sigma_{\bar{x}} (V, ~ \alpha_1 \wedge  ... \wedge \alpha_k) ~ = ~
(V, ~~ D_{\bar{x}}(\alpha_1 \wedge \cdots \wedge \alpha_k).
\]
It follows then that for any $l_i = (\bar{x}_i, \alpha_i) \in
\Gamma(L), i = 1, ..., k$
\beqa
(d_{\scriptscriptstyle L} \theta(V)) (l_1 \wedge \cdots \wedge l_k) & = &
\sum_{i<j} (-1)^{i+j} (\theta(V), ~ \{l_i, l_j\} \wedge l_1 \wedge \cdots
\hat{l}_i \cdots
\hat{l}_j \cdots  \wedge l_k) \\
& = & \sum_{i<j} (-1)^{i+j} (V, ~ \{\alpha_i, \alpha_j\} \wedge \alpha_1
\wedge \cdots
\hat{\alpha}_i \cdots  \hat{\alpha}_j \cdots \wedge  \alpha_k )\\
& & + \sum_{i<j} (-1)^{i+j} (V, ~ (D_{\bar{x}_i}
\alpha_j - D_{\bar{x}_j} \alpha_i) \wedge  \alpha_1 \wedge \cdots
\hat{\alpha}_i \cdots
\hat{\alpha}_j \cdots \wedge \alpha_k )\\
& = & \sum_{i<j} (-1)^{i+j} (V, ~ \{\alpha_i, \alpha_j\} \wedge \alpha_1
\wedge \cdots
\hat{\alpha}_i \cdots  \hat{\alpha}_j \cdots \wedge  \alpha_k )\\
& & + \sum_j (-1)^{j} (V, ~ D_{\bar{x}_j} (\alpha_1 \wedge
\cdots \hat{\alpha}_j \cdots \wedge \alpha_k))\\
& = & \sum_{i<j} (-1)^{i+j} (V, ~ \{\alpha_i, \alpha_j\}\wedge \alpha_1
\wedge \cdots
\hat{\alpha}_i \cdots  \hat{\alpha}_j \cdots \wedge  \alpha_k )\\
& & + \sum_j (-1)^{j+1} \sigma_{\bar{x}_j}(V, ~
\alpha_1 \wedge \cdots \hat{\alpha}_j \cdots \wedge \alpha_k))\\
& = & \sum_{i<j} (-1)^{i+j} (V, ~ \{\alpha_i, \alpha_j\} \wedge \alpha_1
\wedge \cdots
\hat{\alpha}_i \cdots  \hat{\alpha}_j \cdots \wedge  \alpha_k )\\
& & + \sum_j (-1)^{j+1} (-\ps(\alpha_j)) (V, ~ \alpha_1 \wedge
\cdots \hat{\alpha}_j \cdots \wedge \alpha_k))\\
& = & [\pi, ~ V] (\alpha_1 \wedge \cdots \wedge \alpha_k)\\
& = & \theta([\pi, ~ V]) (l_1 \wedge \cdots \wedge l_k).
\eeqa
Therefore
\[
d_{\scriptscriptstyle L} \theta(V) ~ = ~ \theta([\pi, ~ V]).
\]
This finishes the proof of the theorem.
\qed

\begin{exam}
\label{exam_coiso-2}
{\em We continue with Example \ref{exam_coiso-1}. Let $G$ be a connected
Poisson
Lie group, and let $H \subset G$ (denoted by $G_p$ in Example
\ref{exam_coiso-1}) be a connected closed subgroup
of $G$ with Lie algebra ${\frak h}$ such that
\[
\fh^{\perp} ~ = ~ \{\xi \in \fgs: ~~ (\xi, x) = 0, ~~ \forall x \in \fh \}
\subset \fgs
\]
is a Lie subalgebra of $\fgs$.  In this case, there is a unique
Poisson structure on the quotient space $G/H$ such that
the projection from $G$ to $G/H$ is a Poisson map, and the left action
of $G$ on $G/H$ by left translations is a Poisson action.
We have seen in Example \ref{exam_coiso-1} that the Lie algebra
${\frak l}_p$ corresponding to the point $p = eH$ in $G/H$
is the Lie subalgebra $\fh + \fh^{\perp}$ of $\fd =
\fg \bowtie \fgs$. By Theorem \ref{thm_invariant-relative}, we know
that the $G$-invariant Poisson cohomology of $G/H$ is isomorphic to the
relative Lie algebra cohomology of $\fh + \fh^{\perp}$ relative to the
Lie subalgebra $\fh$. We also observed in Example \ref{exam_coiso-1}
that the Lie algebra $\fh + \fh^{\perp}$ is an example of
a double Lie algebra. We first prove the following fact on
relative cohomology in the case of a double Lie algebra.
}
\end{exam}

\begin{lem}
\label{lem_relative-double}
Let ${\frak l}$ be a Lie algebra and let $\fh$ and $\fn$ be Lie subalgebras
of $\fl$ such that
\[
\fl ~ = ~ \fh ~ \oplus ~ \fn
\]
as vector spaces. For $x \in \fh$ and $\xi \in \fn$, write
\[
[x, ~ \xi ] ~ = ~ -\xi \cdot x ~ + ~ x \cdot \xi
\]
where $\xi \cdot x \in \fh$ and $x \cdot \xi \in \fn$. The maps
\begin{eqnarray}
\label{eq-action-1}
& & \fh \times \fn \lrw \fn: ~~ (x, ~ \xi) \Map x \cdot \xi\\
\label{eq_action-2}
& & \fn \times \fh \lrw \fh: ~~ (\xi, ~ x) \Map \xi \cdot x
\end{eqnarray}
define a pair of actions of $\fh$ and $\fn$ on each other. Denote
by $(\wedge \fn^*)^{\fh}$ the subspace of $\wedge \fn^*$ of
$\fh$-invariant vectors with respect to the action of $\fh$ on $\fn^*$
contragradient to the action
of $\fh$ on $\fn$. Let
\[
d_{\fn}: ~ \wedge^{k} \fn^* \lrw \wedge^{k+1} \fn^*
\]
be
the Chevalley-Eilenberg coboundary operator  for the Lie algebra
$\fn$. Then the subspace $(\wedge \fn^*)^{\fh}$ of
$\wedge \fn^*$ is invariant under
$d_{\fn}$ (even though the action of $\fh$ on $\fn$ is not by
Lie algebra derivations).  The cohomology  $H^{\bullet}
((\wedge \fn^*)^{\fh}, ~ d_{\fn})$ of the
cochain subcomplex $((\wedge \fn^*)^{\fh}, ~ d_{\fn})$ is isomorphic
to the relative Lie algebra cohomology of $\fl$ relative to the
Lie subalgebra $\fh$.
\end{lem}

\noindent
{\bf Proof.} By embedding ${\frak n}^*$ into ${\frak l}^*$
via
\[
\fn^* \ni f \Map \bar{f} (x + y) ~ = ~ f(y), ~~~ x \in \fh, ~ y \in \fn
\]
we can identify $\fn^*$ with $\fh^{\perp} \subset \fl^*$ and
$(\wedge \fn^*)^{\fh}$ with $(\wedge \fh^{\perp})^{\fh} \subset \wedge \fl^*$.
Denote again
by $f$ an arbitrary element in $(\wedge \fn^*)^{\fh}$ and by
$\bar{f}$ its image in $(\wedge \fh^{\perp})^{\fh} \subset \wedge \fl^*$, i.e.,
\[
\bar{f} (x_1 + y_1, \cdots, x_k + y_k) ~ = ~ f(y_1, \cdots, y_k).
\]
We need to show that
\[
d_{\fl} \bar{f} ~ = ~ \overline{d_{\fn} f},
\]
where $d_{\fl}: ~ wedge^k \fl^* \rightarrow \wedge^{l+1} \fl^*$ is the
Chevalley-Eilenberg coboundary operator for $\fl$.
Let $l_i = x_i + y_i \in \fh + \fn$ with $x_i \in \fh$ and
$y_i \in \fn, ~ i = 1, ..., k+1$. Then
\beqa
(d_{\fl} \bar{f}) (l_1, ..., l_{k+1}) & = & \sum_{i<j} (-1)^{i+j}
\bar{f}([l_i, ~ l_j], \cdots \hat{l}_i, \cdots, \hat{l}_j, \cdots, l_{k+1})\\
& = & \sum_{i<j} (-1)^{i+j} f([y_i, ~ y_j] + x_i \cdot y_j - x_j \cdot
y_i, \cdots, \hat{y}_i, \cdots, \hat{y}_j, \cdots, y_{k+1} ) \\
& = & \overline{d_{\fn} f}(l_1, \cdots, l_{k+1}) ~ + ~
\sum_j (-1)^j f(x_j \cdot (y_1
\wedge \cdots \wedge \hat{y}_j \wedge \cdots \wedge y_{k+1}))\\
& = & \overline{d_{\fn} f} (l_1, \cdots, l_{k+1}).
\eeqa
In the last step we use that fact that $f$ is $\fh$-invariant.
\qed

\bigskip
The following Theorem now follows from
Lemma \ref{lem_relative-double}.

\begin{thm}
\label{thm_coisotropic-cohomo}
Let $G$ be a connected Poisson Lie group, and let $H \subset G$
be a connected closed subgroup
of $G$ with Lie algebra ${\frak h}$ such that
\[
\fh^{\perp} ~ = ~ \{\xi \in \fgs: ~~ (\xi, x) = 0, ~~ \forall x \in \fh \}
\subset \fgs
\]
is a Lie subalgebra of $\fgs$. Denote by $(\wedge (\fh^{\perp})^*)^{\fh}$
the subspace in $\wedge (\fh^{\perp})^*$ of $\fh$-invariant
vectors with respect to the adjoint action of $\fh$ on
$(\fh^{\perp})^* \cong \fg/\fh$. Let
\[
d_{\fh^{\perp}}: ~~\wedge^k (\fh^{\perp})^* \lrw \wedge^{k+1}
(\fh^{\perp})^*
\]
be the Chevalley-Eilenberg coboundary operator for the Lie algebra
$\fh^{\perp}$. Then $(\wedge (\fh^{\perp})^*)^{\fh}$ is
invariant under $d_{\fh^{\perp}}$, and the the cohomology of
$((\wedge (\fh^{\perp})^*)^{\fh}, ~ d_{\fh^{\perp}})$
is isomorphic to the $G$-invariant Poisson cohomology
of the Poisson homogeneous space $G/H$.
\end{thm}

\bigskip
We now apply Theorem \ref{thm_coisotropic-cohomo} to the
Bruhat Poisson structure on a generalized flag manifold \cite{lu-we:poi}.

Let $G$ be a connected finite dimensional semisimple complex Lie group
and let $G = KAN$ be an Iwasawa decomposition of $G$ as
a real semisimple Lie group. Then there are natural Poisson
structures on $K$ and on $AN$ making them into Poisson
Lie groups. Moreover, as Poisson Lie groups, they are
dual to each other, and the group $G$ is the double of $K$
and $AN$. On the Lie algebra level, we have
the Iwasawa decomposition $\fg = {\frak k} + {\frak a} +
{\frak n}$. The imaginary part of the (complex-valued)
Killing form of $\fg$ gives a nondegenerate pairing
between ${\frak k}$ and ${\frak a} + {\frak n}$. In
other words, the pair $({\frak k}, ~ {\frak a} + {\frak n})$
form a Lie bialgebra and $\fg$ is its double Lie algebra.

Let $T$ be the maximal torus in $K$, and let ${\frak t}$
be its Lie algebra. The subspace
\[
{\frak t}^{\perp} ~ = ~ \{\xi \in {\frak k}^*: ~ (\xi, x) = 0, ~~
\forall x \in {\frak t} \}
\]
of ${\frak k}^* \cong {\frak a} + {\frak n}$ is clearly
${\frak n}$ which is an ideal in ${\frak a} + {\frak n}$.
Thus there is a unique Poisson structure on $K/T$ such
that the projection from $K$ to $K/T$
is a Poisson map. Moreover, the left action of $K$ on $K/T$ is
a Poisson action, so $K/T$ is a Poisson homogeneous $K$-space.

Symplectic leaves in $K/T$ for this Poisson structure are
known to be exactly the Bruhat cells in $K/T$ (thus the name
Bruhat Poisson structure).
This can be
seen directly from Theorem \ref{thm_leaves}, for in this case,
the Lie group $L_p$, for $p = eT \in K/T$, is $TN$, and the
right action of it on $K$  described in Theorem \ref{thm_leaves}
is just the restriction to $TN$ of the right action of $G$ on
$K \cong (AN)\backslash G$ by right translations. The
projection of these orbits from $K$ to $K/T$ are exactly
the Bruhat cells in $K/T$.

Consider now the $K$-invariant Poisson cohomology space of the
Bruhat Poisson structure with complex coefficients. By
Theorem \ref{thm_coisotropic-cohomo}, we get the following

\begin{thm}
\label{thm_bruhat}
The (complex-valued) $K$-invariant Poisson cohomology
$H_{\pi, {\scriptscriptstyle K}}(K/T)$
of the Bruhat Poisson structure on $K/T$
is isomorphic to the space $End_{\fh}(H(\fn))$,
where $\fh$ is the complexification of ${\frak t}$, and
$\fn$ here is considered as a complex Lie algebra.
Thus the dimension of $H_{\pi, {\scriptscriptstyle K}}^{k}(K/T)$
is $0$ if $k$ is odd and is equal to the number of Weyl group element
of length $k/2$ if $k$ is even. In particular, it is isomorphic to the
(complex-valued) de Rham cohomology of $K/T$.
\end{thm}

\bigskip
\noindent
{\bf Proof.} From Theorem \ref{thm_coisotropic-cohomo} we know that
\[
H_{\pi, {\scriptscriptstyle K}}(K/T) ~ \cong ~ H( (\wedge
({\frak n} + {\frak n}_{-})^*)^{\fh}, ~ d_{\fn}) ~ \cong ~ End_{\fh}(H(\fn)).
\]
By Kostant's theorem \cite{kos:cohomology},
there is exactly one cohomology class for the Lie algebra ${\frak n}$
corresponding to each element in the Weyl group $W$ of $K$. Choose
root vectors $\{E_{\alpha}, E_{-\alpha}: \alpha > 0\}$.
For each $w \in W$, set
\[
\Phi_w ~ :=~ \{ \beta > 0: ~ w^{-1} \beta > 0 \} ~ = ~ \{\beta_1, \cdots,
\beta_k\}.
\]
Then the (complex) $K$-invariant $2k$-vector field
\[
E_{-\beta_1} \wedge \cdots \wedge E_{-\beta_k} \wedge E_{\beta_1}
\wedge \cdots \wedge E_{\beta_k}
\]
on $K/T$ is a representative of a cohomology class of
the $K$-invariant
Poisson cohomology of the Bruhat Poisson structure on $K/T$.
The fact that the space $End_{\fh}(H(\fn))$ is isomorphic to
the de Rahm cohomology of $K/T$ is another theorem of Kostant's
\cite{kos:schubert}.

\section{The Lie groupoid of $A$}
\label{sec_groupoid}

We first recall the definition of a Lie groupoid. Details can be found in
\cite{mk:book}.

A groupoid over a set $P$ is a set $M$, together with

 (1)  surjections $s, t: M \rightarrow P$ (called
the source and target maps respectively);

 (2)  $\mu: M * M  \rightarrow M $ (called the multiplication), where
\[
M * M = \{(m_2, m_1) \in M \times M: ~ s(m_2) = t(m_1)\};
\]
each pair $(m_2, m_1)$ in $M*M$ is said to be ``composable'';

 (3) an injection $\epsilon: P \rightarrow M$ (identities);

 (4) $\iota: M \rightarrow M$ (inversion).

These maps must satisfy

 (1) associative law: $\mu (\mu (m_3, m_2), m_1) ~ = ~ \mu(m_3, \mu(
m_2, m_1))$ (if one
side is defined, so is the other);

 (2) identities: for each $m \in M, ~ (\epsilon(t(m)), m) \in M*M, ~
(m, \epsilon (s(m)) \in M*M$, and $\mu(\epsilon(t(m)), m)
 = \mu (m, \epsilon (s(m)) = m$;

 (3) inverses: for each $m \in M, (m, \iota(m)) \in M*M,~  (\iota (m),
 m) \in M*M$, and $\mu(m, \iota(m)) = \epsilon (t(m))$ and
$\mu(\iota (m), m) = \epsilon (s(m)).$

A Lie groupoid (or differential groupoid) $M$ over a manifold
$P$ is a groupoid with a differential structure such that
(1) $s$ and $t$ are differentiable submersions (this implies that
$M*M$ is a submanifold of $M \times M$), and
(2)  $\mu, \epsilon$ and $\iota$ are differentiable maps.

Given a Lie groupoid $M$ over $P$, the normal bundle
of $P$ in $M$ has a Lie algebroid structure over
$P$ whose sections can be identified with left invariant
vector fields on $M$. It is called the Lie algebroid of $M$.

A left action of a Lie groupoid  $(M, P, s,t)$ on a smooth submersion
$f: Q \rightarrow P$
is a map
\[
M *_f Q ~ = ~ \{(m, q): ~ s(m) = f(q)\} \lrw Q: ~~ (m, ~ q) \Map mq
\]
such that

(1) $f(mq) = t(m)$

(2) $m_2 (m_1 q) ~ = ~ (m_2m_1) q$

(3) $\epsilon(f(q)) q = q$

for all $m, m_1, m_2 \in M$ and $q \in Q$ which are suitably compatible.

Given such an action, one can construct the transformation groupoid
on $M *_f Q$ with
base $Q$ by defining
\beqa
& & s'(m,~q) = q, ~~~~ t'(m,~q) = mq\\
& & (m_2, ~ m_1 q) (m_1, ~ q) ~ = ~ (m_2 m_1, ~ q).
\eeqa
The map
\[
M *_f Q \rightarrow M: ~ (m, ~ q) \Map m
\]
is a morphism of Lie groupoid over $f: Q \rightarrow P$.

A right action of a groupoid and the corresponding
transformation groupoid are similarly defined.

\begin{exam}
\label{exam_group}
{\em
A Lie group is a Lie groupoid over a one point space. A Lie
group action $G \times P \rightarrow P$ is an example
of a Lie groupoid action. The Lie algebroid of the
corresponding transformation groupoid is the
transformation Lie algebroid discussed in Section \ref{sec_lie-algebroids}.
}
\end{exam}

\begin{exam}
\label{exam_symplectic-groupoid}
{\em
Let $P$ be a Poisson manifold. We have seen in Section \ref{sec_lie-algebroids}
that the cotangent bundle $T^*P$ of $P$ has a Lie algebroid structure. If there
is a Lie groupoid $(N, \sn, \tn)$ over $P$ whose Lie
algebroid is $T^*P$ and whose $\sn$-fibers are simply-connected, we say that
$P$ is integrable as a Poisson manifold. It turns out that
there is always a symplectic structure on $N$ making it into
a symplectic groupoid \cite{cdw:gpoid}. We call $N$ the
symplectic groupoid  of $P$.
}
\end{exam}

\bigskip
Assume now that $G$ is a Poisson Lie group, $P$ is an integrable Poisson
manifold with symplectic groupoid $(N, P, \sn, \tn)$,
and $\sigma: G \times P \rightarrow P$ is
a Poisson action. Denote by $(M, \sm, \tm)$ the
transformation Lie groupoid over $P$ defined by $\sigma$.
The purpose of this section is to describe a Lie groupoid whose
Lie algebroid is the Lie algebroid $A = (P \times \fg) \bowtie
T^*P$ given in Theorem \ref{thm_A}. Not surprisingly, this
groupoid will be built out of $M$ and $N$ and a pair of actions of $M$
and $N$ on each other.

\bigskip
Recall that the map
\[
\phi: ~ T^*P \lrw \fgs: ~~ (\phi(\alpha_p), ~ x) ~ = ~ (\alpha_p, ~ \sigma_x
(p))
\]
is a Lie algebroid homomorphism. Assume that it can be integrated
to a Lie groupoid homomorphism
\[
\varphi: ~~ N \lrw G^*,
\]
where $G^*$ is the dual Poisson Lie group of $G$
(see \cite{xu:poisson-groupoid}).
The map $\varphi$ is aslo a Poisson map. As a such, it
induces a left action of $G$ on $N$ \cite{lu:mom} which we will denote by
$(g, n) \Map gn$. The map $\tp$ and the action of $G$ on $N$
have the following properties: for any $g \in G, ~ n, n_1, n_2 \in N$,
\beqa
& & \tn(gn) ~ = ~ g \tn(n)\\
& & \sn(gn) ~ = ~ g^{\tp(n)} \sn(n)\\
& & \en(gn) ~ = ~ g \en(n)\\
& & g(n_2n_1) ~ = ~ (g n_2) (g^{\tp(n_2)} n_1)\\
& & \tp(\en(p)) ~ = ~ 1 \in G^*\\
& & \tp (n_2n_1) = \tp(n_2) \tp(n_1)\\
& & \tp(gn) ~ = ~ ^{g} \! \tp(n).
\eeqa
Here $^{g}\! u$ and $g^u$, for $g \in G$ and $u \in G^*$, denote
respectively the left action of $G$ on $G^*$ and the right
action of $G^*$ on $G$ defined by the decomposition
\[
gu ~ = ~ ^{g}\!u ~g^u \in D
\]
in the group $D$ with Lie algebra $\fd$ (see Section \ref{sec_poi}).

\bigskip
We can now describe the pair of actions of the groupoids $M$ and $N$
on each other.

The left action of $M$ on $\tn: N \rightarrow P$
is given by
\begin{equation}
\label{eq_M-on-N}
M *_{\tn} N \lrw N: ~~ ((g, ~ \tn(n)), ~ n) \Map gn.
\end{equation}
The right action of $N$ on $M$ is given by
\begin{equation}
\label{eq_N-on-M}
M *_{\sm} N ~ = ~ M *_{\tn} N: ~~ ((g, ~ \tn(n)), ~ n) \Map
(g^{\tp(n)}, ~ \sn(n)).
\end{equation}
For notational simplicity, we denote the manifold $M *_{\tn} N
= M *_{\sm} N$ by $M * N$ and denote the two actions respectively by
\beqa
& & M * N \lrw N: ~~ (m, ~ n) \Map ^{m}\!n\\
& & M * N \lrw M: ~~ (m, ~ n) \Map m^{n}.
\eeqa
This pair of actions are compatible in the following sense:
\beqa
& & (1) ~ ~\sn(^{m}\! n) ~ = ~ \tm (m^n), ~~~ \forall (m, n) \in M*N\\
& & (2) ~ ~^{m}\! (\en(\sm(m))) ~ = ~ \en (\tm(m)), ~~~ \forall m \in M \\
& & (3) ~~ (\mem (\tn(n)))^n ~ = ~ \mem(\sn(n)), ~~~\forall n \in N\\
& & (4) ~~ ^{m}\! (n_2n_1) ~ = ~ (^{m}\! n_2) (^{m^{n_2}}\! n_1), ~~~ \forall
(m, n_2) \in M*N, ~ (n_2, n_1) \in N*N\\
& & (5) ~~ (m_2m_1)^n ~ = ~ (m_{2}^{^{m_1}\!n}) (m_{1}^{n}), ~~~ \forall (m_2,
m_1) \in M*M.
\eeqa

According to \cite{mk:double}, two groupoids $M$ and $N$ over the
same base $P$  with a pair of actions of them on each other satisfying
conditions (1)-(5) as above are called a pair of groupoids with
an interaction or a ``matched pair" of groupoids.
Given such data, one can construct a third groupoid
on the space $N*M = \{(n, m): ~ \sn(n) = \tm(m)\}$ over $P$
by defining
\beqa
& & s(n, m) ~ = ~ \sm(m), ~~~ t(n, m) ~ = ~ \tn (n)\\
& & (n_1, m_1) (n_2, m_2) ~ = ~ (n_1~(^{m_1}\! n_2) , ~~ (m_{1}^{n_2})~ m_2 ).
\eeqa
The natural maps
\beqa
& & M \lrw N*M: ~~ m \Map (\en(\tm(m)), ~ m)\\
& & N \lrw N*M: ~~ n \Map (n, ~ \mem(\sn(n))
\eeqa
embed $M$ and $N$ into $N*M$ as subgroupoids.

\begin{thm}
\label{thm_groupoid-for-A}
Let $\sigma: G \times P \rightarrow P$ be a Poisson action of
a Poisson Lie group $G$ on an integrable Poisson manifold $P$.
Let $N$ be the symplectic groupoid of $P$.
Let $N*M$ be the Lie groupoid over $P$ as described above.
Then the Lie algebroid of $N*M$ is the Lie algebroid
$A = (P \times \fg) \bowtie T^*P$ as given in
Theorem \ref{thm_A}.
\end{thm}

We will skip the proof of this theorem. The main part of the proof
consists of showing
that the linearizations of the two actions of $M$ and $N$
on each other at their identity sections
are exactly the Lie algebroid representations of $P \times \fg$ and
$T^*P$ on each other which were used in the construction
of the Lie algebroid $A = (P \times \fg) \bowtie T^*P$.

\begin{rem}
\label{rem_fiber}
{\em
For $p \in P$, denote by $(N*M)_p$ the intersection of the
$s$-fiber and the $t$-fiber in $N*M$ over $p$. It can be
identified with the space
\[
(N*M)_p ~ = ~ \{(n, ~ g): ~ n \in N , ~ g \in G: ~ \tn(n) = p, ~
\sn(n) = gp \}.
\]
We know from general groupoid theory that it has a group
structure:
\[
(n_2, ~ g_2) (n_1, ~ g_1) ~ = ~ (n_2 (g_2 n_1), ~~
g_{2}^{\varphi(n_1)} ~ g_1), ~~~~
(n_2, ~ g_2), ~ (n_1, ~ g_1) \in (N*M)_p.
\]
The Lie algebra of $(N*M)_p$ is the Lie algebra
${\frak l}_p$ introduced in Section \ref{sec_A}.
}
\end{rem}

\begin{rem}
\label{rem_double}
{\em
We remark that the space $N*M$ is an example of what is called
a vacant double groupoid in \cite{mk:double}. On the space $N *M$,
there is also a groupoid structure over $M$ and a groupoid
structure over $N$, namely the transformation groupoids of the actions
of $M$ and $N$ on each other. These two groupoid structures
are compatible, making $N*M$ into a double groupoid. The groupoid structure
on $P$ we just described is called the diagonal groupoid of the
double groupoid $M*N$. See \cite{mk:double} for more
details. For a detailed treatment of ``matched pairs" of Lie algebroids,
see \cite{mk:thesis}.
}
\end{rem}

\end{document}